\documentclass[preprint]{amsart}

\usepackage{graphicx,subfigure}
\usepackage{epstopdf}
\usepackage{amsmath,amsthm,amssymb}
\usepackage{rotating}
\usepackage{morefloats}
\usepackage{comment}
\usepackage{amsmath}

\newtheorem{prp}{Proposition}
\newtheorem{thm}{Theorem}
\newtheorem{crl}{Corollary}

\begin{document}
\title{Information geometry of tempered stable processes}
\author{Jaehyung Choi}
\address{}
\email{jj.jaehyung.choi@gmail.com}

\begin{abstract}
	We derive the information geometry of tempered stable processes. We first compute the $\alpha$-divergence between two tempered stable processes. From the divergence, we obtain the Fisher information matrix and the $\alpha$-connection on the statistical manifold. The generalized, classical, and rapidly-decreasing tempered stable geometries are dually flat, with potential functions and canonical divergences expressed in terms of the tempering parameters. Moreover, their $\alpha$-curvature tensors vanish for every $\alpha$. The geometric results are also applied to bias reduction in maximum likelihood estimation and Bayesian predictive priors.
\end{abstract}

\maketitle
\section{Introduction}
\label{sec_intro}
	Information geometry is an interdisciplinary study across differential geometry, probability theory, statistics, and information theory. With the mathematically elegant structure of Riemannian geometry, information geometry provides new interpretations on statistics and probability theory \cite{efron1975defining, efron1978geometry,firth1993bias,kosmidis2009bias,kosmidis2010generic}. Beyond its mathematical elaboration, the differential geometry of statistical manifolds also offers many practical applications to various probabilistic models \cite{komaki2006shrinkage,tanaka2008superharmonic,tanaka2018superharmonic, choi2015kahlerian, choi2015geometric, oda2021shrinkage}.
	
	Time series analysis and signal processing have been empowered by such applications of information geometry. Since the derivation for the information geometry of time series models such as autoregressive and moving average (ARMA) models \cite{ravishanker1990differential} and fractionally-integrated ARMA (ARFIMA) models \cite{ravishanker2001differential}, there have been several approaches to theoretical advances for the statistical manifolds of time series models and signal processing filters by introducing symplectic geometry and K\"ahler geometry \cite{barndorff1997statistics, barbaresco2006information, barbaresco2012information, zhang2013symplectic, barbaresco2014koszul, choi2015kahlerian,choi2021kahlerian}. Moreover, the mathematical sophistication and enhancements produce practical advantages such as Bayesian predictive priors outperforming the Jeffreys prior and their easier derivation \cite{komaki2006shrinkage, tanaka2008superharmonic, tanaka2018superharmonic, choi2015kahlerian, choi2015geometric, oda2021shrinkage}.
	
	Applications of information geometry are not restricted to time series analysis and signal processing. Information geometry also plays a pragmatic role in statistical inference and parameter estimation for probabilistic models. For instance, the concept of statistical curvature for exponential families was introduced based on differential geometry \cite{efron1975defining, efron1978geometry}. Additionally, the differential-geometric approach allows for bias reduction in parameter estimation \cite{firth1993bias,kosmidis2009bias,kosmidis2010generic}. Similar to time series analysis and signal processing, we also can use information geometry to find outperforming Bayesian predictive priors for probability distributions \cite{komaki2006shrinkage}. 
	
	However, not all probability distributions and statistical models are eligible for these benefits and advantages of information geometry because we may have technical difficulties and limitations in exploring the information geometry of several probability distributions and statistical models. For example, some distributions may lack closed-form expressions for their probability density functions. Even when an explicit expression for the probability density function is available, the Fisher information matrix could be infinite or difficult to compute, and deriving the metric tensor of information geometry can still be technically complex.
	
	The classical tempered stable (CTS) distribution, suggested by Rosinski \cite{rosinski2007tempering} and also known as CGMY distribution \cite{carr2002fine}, is one such probability distribution for which investigating the information geometry is difficult. Although the distribution has been applied in various financial contexts such as financial time series modeling using the ARMA-GARCH-CTS model \cite{kim2009computing, kim2010tempered, kim2011time}, as well as in portfolio management and options pricing \cite{tsuchida2012mean, beck2013empirical, choi2015reward, georgiev2015periodic, anand2017equity, kim2023deep, choi2024diversified}, the closed-form expression of its probability density function does not exist. This lack of the closed-form probability density function presents the fundamental obstacle in deriving the information geometry of CTS processes. The same challenge applies not only to the general version of CTS processes, called generalized tempered stable (GTS) processes \cite{rachev2011financial}, but also to other tempered stable processes such as rapidly-decreasing tempered stable (RDTS) processes \cite{kim2010tempered} and a broad class of L\'evy processes \cite{rachev2011financial}.
	
	Ongoing efforts to explore the information geometry of tempered stable processes can be found in the work of Kim and Lee \cite{kim2007relative} where they derived the relative entropy between two equivalent martingale CTS processes from the Radon--Nikodym derivative between the corresponding L\'evy measures. However, their study does not provide the information geometry of CTS processes. Furthermore, the statistical manifolds of both CTS processes and GTS processes remain unexplored, and even the relative entropy of GTS processes is unknown.
	
	In this paper, we explore the information geometry of tempered stable processes. Generalizing the result found by Kim and Lee \cite{kim2007relative} not only from relative entropy to $\alpha$-divergence but also from CTS processes to other tempered stable processes, we calculate the Fisher information matrices of tempered stable processes and the $\alpha$-connections of the tempered stable process geometry. Additionally, we also derive the information geometry of RDTS processes. We further study dual flatness, potential functions, canonical divergences, and curvature of these geometries. The geometric results are also applied to bias reduction and Bayesian predictive priors.
	
	The structure of this paper is as follows. In the next section, we provide the fundamentals of tempered stable processes, and re-visit Kim and Lee's work on the relative entropy of CTS processes \cite{kim2007relative}. Section~\ref{sec_ig} provides the main results of this paper such as the information geometry of various tempered stable processes. In Section~\ref{sec_apps}, statistical applications for the information geometry of GTS processes are covered. Finally, we conclude the paper.
	
\section{Tempered stable processes}
	In this section, we provide the basic concepts of various tempered stable processes such as CTS processes, GTS processes, and RDTS processes. We also review the prior work on the relative entropy of CTS processes. These materials will be the foundation for deriving our main results in the next section.

	A tempered stable distribution is an infinitely divisible distribution with the L\'evy triplet $(0,\nu,\gamma)$, where $\nu$ is a tempered stable L\'evy measure \cite{rachev2011financial}. The tempered stable L\'evy measure is obtained by tempering a stable-type L\'evy measure \cite{rosinski2007tempering} as
	\begin{align}
	\label{levy_measure_ts}
		\nu(dx;\boldsymbol{\xi})=t(x;\boldsymbol{\xi})\nu_{\textrm{stable}}(dx;\boldsymbol{\xi})=\lambda(x;\boldsymbol{\xi})dx,
	\end{align}
	where $t(x;\boldsymbol{\xi})$ is a tempering function. The stable-type L\'evy measure used in this paper is represented with
	\begin{align}
	\label{levy_measure_stable}
		\nu_{\textrm{stable}}(dx;\boldsymbol{\xi})=\Big(\frac{C_+}{x^{a_++1}} 1_{x>0}(x)+\frac{C_-}{|x|^{a_-+1}} 1_{x<0}(x)\Big)dx,
	\end{align}
	where $C_+, C_-$ are positive, $a_+, a_-\in(0,2)\setminus \{1\}$, and $1_A(x)$ is the indicator function such that $1_A(x)=1$ when $x\in A$ and 0 otherwise. 
	
	The characteristic function of the tempered stable distribution with the L\'evy triplet $(0,\nu,\gamma)$ is obtained by the L\'evy--Khintchine formula \cite{rachev2011financial} as
	\begin{align}
	\label{charac_ftn_ts}
		\phi(u;\gamma,\boldsymbol{\xi})=\exp\Big(iu\gamma+\int_{-\infty}^{\infty}({\rm e}^{iux}-1-iux 1_{|x|\leq1})\nu(dx;\boldsymbol{\xi})\Big).
	\end{align}
	
	A tempered stable process is a L\'evy process whose unit-time distribution is a tempered stable distribution. Equivalently, its L\'evy measure has the form in Eq.~(\ref{levy_measure_ts}). For a tempered stable process $(X_t,\mathbb{P})_{t\in[0,T]}$ with $X_0=0$, the characteristic function of $X_t$ is
	\begin{align}
		\mathbb{E}[{\rm e}^{iuX_t}]=\exp\Big(t\big(iu\gamma+\int_{-\infty}^{\infty}({\rm e}^{iux}-1-iux 1_{|x|\leq1})\nu(dx;\boldsymbol{\xi})\big)\Big).
	\end{align}
	
	More details on tempered stable distributions and processes can be found in Rachev et al. \cite{rachev2011financial} and references therein.
	
\subsection{CTS processes and GTS processes}
	Since Rosinski \cite{rosinski2007tempering} introduced the concept of tempering stable distributions, which led to the development of the classical tempered stable (CTS) distribution, widespread applications have been found across various fields. In particular, the CTS distribution, also known as the CGMY distribution \cite{carr2002fine}, is well-suited for modeling heavy tails and skewness. As a result, it has been extensively applied in finance, including time series modeling through the ARMA-GARCH-CTS model \cite{kim2009computing, kim2010tempered, kim2011time}, as well as in portfolio management and options pricing \cite{tsuchida2012mean, beck2013empirical, choi2015reward, georgiev2015periodic, anand2017equity, kim2023deep, choi2024diversified}.
	
	 A CTS distribution is characterized by five parameters, $\boldsymbol{\xi}=(a, C,\lambda_+,\lambda_-,m)$: $m$ is the location parameter; $a$ is the tail index; $C$ is the scale parameter; $\lambda_+$ and $\lambda_-$ are the decay rates of the upper and lower tails, respectively. $C, \lambda_+, \lambda_-$ are positive, $a\in(0,2)\setminus \{1\}$, and $m\in\mathbb{R}$.
	 
	By imposing the following CTS condition on tempered stable distributions in Eq.~(\ref{levy_measure_ts}):
	\begin{align}
	\label{cts_condition}
		\left\{ 
		\begin{array}{ll}
			a_+=a_-=a\\
			C_+=C_-=C
		\end{array}
		\right.,
	\end{align}
	and using the CTS tempering function of
	\begin{align}
	\label{temp_fns_cts}
	 	t_{\textrm{CTS}}(x;\boldsymbol{\xi})={\rm e}^{-\lambda_{+}x} 1_{x>0}(x)+{\rm e}^{-\lambda_{-} |x|} 1_{x<0}(x),
	\end{align}
	the L\'evy measure of the CTS distribution from Eq.~(\ref{levy_measure_ts}) \cite{rachev2011financial} is given by
	\begin{align}
	\label{levy_measure_cts}
		\nu(dx;\boldsymbol{\xi})=C\Big(\frac{{\rm e}^{-\lambda_{+}x}}{x^{a+1}} 1_{x>0}(x)+\frac{{\rm e}^{-\lambda_{-} |x|}}{|x|^{a+1}} 1_{x<0}(x)\Big)dx.
	\end{align}
	
	Based on Eq.~(\ref{charac_ftn_ts}) and Eq.~(\ref{levy_measure_cts}), the characteristic function of the CTS distribution is represented with
	\begin{align}
	\label{charac_ftn_cts}
	\begin{split}
		\phi(u;\boldsymbol{\xi})=&\exp\Big(ium-iuC\Gamma(1-a)(\lambda_+^{a-1}-\lambda_-^{a-1})\\
		&+C\Gamma(-a)\big(\big((\lambda_+-iu)^a-\lambda_+^a\big)+\big((\lambda_-+iu)^a-\lambda_-^a\big)\big)\Big),
	\end{split}
	\end{align}
	where $\Gamma$ is the Gamma function. Meanwhile, a closed-form expression for the probability density function of CTS distributions does not exist.
		
	Instead of controlling the upper and lower tails by the same tail index and the identical scale parameter, we can introduce additional degrees of freedom to both parameters. The upper and lower tails can be modeled by separate tail indexes and scale parameters. This distribution with these extra parameters is called a GTS distribution \cite{rachev2011financial}. 
	
	We have seven parameters for a GTS distribution, $\boldsymbol{\xi}=(a_+, a_-, C_+, C_-,\lambda_+,\lambda_-,m)$: $m$ is the location parameter; $a_+$ and $a_-$ are the tail indexes; $C_+$ and $C_-$ are the scale parameters; $\lambda_+$ and $\lambda_-$ are the decay rates of the upper and lower tails, respectively. $C_+, C_-, \lambda_+,\lambda_-$ are positive, $a_+, a_-\in(0,2)\setminus \{1\}$, and $m\in\mathbb{R}$.
	
	By using the tempering function of Eq.~(\ref{temp_fns_cts}) for a GTS distribution, i.e.,
	\begin{align}
	\label{temp_fns_gts}
	 	t_{\textrm{GTS}}(x;\boldsymbol{\xi})={\rm e}^{-\lambda_{+}x} 1_{x>0}(x)+{\rm e}^{-\lambda_{-} |x|} 1_{x<0}(x),
	\end{align}
	in Eq.~(\ref{levy_measure_ts}), the L\'evy measure of the GTS distribution \cite{rachev2011financial} is given by
	\begin{align}
	\label{levy_measure_gts}
		\nu(dx;\boldsymbol{\xi})=\Big(\frac{C_+{\rm e}^{-\lambda_+ x}}{x^{a_++1}} 1_{x>0}(x)+\frac{C_-{\rm e}^{-\lambda_- |x|}}{|x|^{a_-+1}} 1_{x<0}(x)\Big)dx.
	\end{align}
	
	From Eq.~(\ref{charac_ftn_ts}) and Eq.~(\ref{levy_measure_gts}), the characteristic function of the GTS distribution \cite{rachev2011financial} is represented with
	\begin{align}
	\label{charac_ftn_gts}
	\begin{split}
		\phi(u;\boldsymbol{\xi})=&\exp\big(ium-iu(C_+\Gamma(1-a_+)\lambda_+^{a_+-1}-C_-\Gamma(1-a_-)\lambda_-^{a_--1})\\
		&+C_+\Gamma(-a_+)\big((\lambda_+-iu)^{a_{+}}-\lambda_+^{a_{+}}\big)+C_-\Gamma(-a_-) \big((\lambda_-+iu)^{a_{-}}-\lambda_-^{a_{-}}\big)\big),
	\end{split}
	\end{align}
	where $\Gamma$ is the Gamma function. Similar to CTS distributions, a closed-form expression for the probability density function of GTS distributions does not exist.
		
	It is straightforward to verify that the CTS distribution is a special case of the GTS distribution. As described above, Eq.~(\ref{cts_condition}) implies that the upper and lower tails of the distribution are controlled by the same tail index, $a$, and the same scale parameter, $C$, but different decay rates, $\lambda_{+}$ and $\lambda_{-}$, respectively. By imposing the CTS condition of Eq.~(\ref{cts_condition}), it is obvious to derive Eq.~(\ref{levy_measure_cts}) and Eq.~(\ref{charac_ftn_cts}) from Eq.~(\ref{levy_measure_gts}) and Eq.~(\ref{charac_ftn_gts}), respectively.
	
	Based on the definition above, a GTS process with parameters $\boldsymbol{\xi}=(a_+,a_-,C_+,C_-,\lambda_+,\lambda_-,m)$ is a L\'evy process with the L\'evy measure in Eq.~(\ref{levy_measure_gts}) and the drift term determined by $m$. Similarly, a CTS process with parameters $\boldsymbol{\xi}=(a,C,\lambda_+,\lambda_-,m)$ is a L\'evy process with the L\'evy measure in Eq.~(\ref{levy_measure_cts}) and the drift term determined by $m$.

	 The next step is to consider the Radon--Nikodym derivative between the L\'evy measures of two GTS processes. Let $(X_t, \mathbb{P})_{t\in[0,T]}$ and $(X_t, \mathbb{Q})_{t\in[0,T]}$ be GTS processes with parameters $\boldsymbol{\xi}=(a_+, a_-, C_+, C_-, \lambda_+, \lambda_-, m)$ and $\tilde{\boldsymbol{\xi}}=(\tilde{a}_+, \tilde{a}_-, \tilde{C}_+,\tilde{C}_-, \tilde{\lambda}_{+}, \tilde{\lambda}_{-}, \tilde{m})$, respectively. We impose the equivalent martingale measure condition, which is based on the L\'evy measures of Eq.~(\ref{levy_measure_gts}) \cite{kim2007relative}.
	
	The equivalent martingale measure (EMM) condition for GTS processes \cite{rachev2011financial} is as follows:
	\begin{align}
	\label{emm_gts}
	\left\{ 
		\begin{array}{ll}
			C_+=\tilde{C}_+\\
			C_-=\tilde{C}_-\\
			a_+=\tilde{a}_+\\
			a_-=\tilde{a}_-\\
			m=\tilde{m}+C_+\Gamma(1-a_+)(\lambda_+^{a_+-1}-\tilde{\lambda}_{+}^{a_+-1})-C_-\Gamma(1-a_-)(\lambda_-^{a_--1}-\tilde{\lambda}_{-}^{a_--1})
		\end{array}
		\right.,
	\end{align}
	and by using Eq.~(\ref{cts_condition}), the equivalent martingale measure condition for CTS processes \cite{rachev2011financial} is given by
	\begin{align}
	\label{emm_cts}
	\left\{ 
		\begin{array}{ll}
			C=\tilde{C}\\ 
			a=\tilde{a}\\
			m=\tilde{m}+C\Gamma(1-a)\big((\lambda_+^{a-1}-\tilde{\lambda}_{+}^{a-1})-(\lambda_-^{a-1}-\tilde{\lambda}_{-}^{a-1})\big)
		\end{array}
		\right..
	\end{align}
	
	With the equivalent martingale measure conditions above, the corresponding Radon--Nikodym derivative for CTS processes and GTS processes \cite{rachev2011financial} is obtained from Eq.~(\ref{levy_measure_cts}) and Eq.~(\ref{levy_measure_gts}), respectively:
	\begin{align}
	\label{radon_nikodym_gts}
		\frac{d\nu^\mathbb{P}}{d\nu^\mathbb{Q}}={\rm e}^{-(\lambda_+-\tilde{\lambda}_+)x} 1_{x>0}(x)+{\rm e}^{-(\lambda_{-}-\tilde{\lambda}_{-})|x|} 1_{x<0}(x),
	\end{align}
	and it is noteworthy that the expression is the same for CTS processes and GTS processes.
	
	For further calculation, it is also convenient to introduce the logarithmic Radon--Nikodym derivative:
	\begin{align}
	\label{log_radon_nikodym}
		\psi(x;\boldsymbol{\xi})=\log{\Bigg(\frac{d\nu^\mathbb{P}}{d\nu^\mathbb{Q}}\Bigg)}.
	\end{align}
	The logarithmic Radon--Nikodym derivative for two equivalent martingale GTS processes is expressed with
	\begin{align}
	\label{log_radon_nikodym_gts}
		\psi(x;\boldsymbol{\xi})=-(\lambda_+-\tilde{\lambda}_+)x 1_{x>0}(x)-(\lambda_{-}-\tilde{\lambda}_{-})|x| 1_{x<0}(x),
	\end{align}
	and it is obvious that the logarithmic Radon--Nikodym derivative for two equivalent martingale CTS processes is also identical to Eq.~(\ref{log_radon_nikodym_gts}).
	
	The Kullback--Leibler divergence, also known as relative entropy, between two tempered stable processes, $(X_t, \mathbb{P})_{t\in[0,T]}$ and $(X_t, \mathbb{Q})_{t\in[0,T]}$, is expressed in terms of the logarithmic Radon--Nikodym derivative in Eq.~(\ref{log_radon_nikodym}) \cite{cont2004nonparametric, kim2007relative}:
	\begin{align}
	\label{kld_measure}
		KL(\mathbb{P}||\mathbb{Q})=T\int (\psi(x){\rm e}^{\psi(x)}-{\rm e}^{\psi(x)}+1)\nu^\mathbb{Q}(dx).
	\end{align}
	
	It is noteworthy that the Kullback--Leibler divergence between two CTS processes can be found from Theorem~2 in Kim and Lee \cite{kim2007relative} based on Eq.~(\ref{kld_measure}).
	\begin{prp}[Kim and Lee (2007)]
		Let $(X_t, \mathbb{P})_{t\in[0,T]}$ and $(X_t, \mathbb{Q})_{t\in[0,T]}$ be CTS processes with parameters $(a, C, \lambda_+, \lambda_-, m)$ and $(\tilde{a}, \tilde{C}, \tilde{\lambda}_{+}, \tilde{\lambda}_{-}, \tilde{m})$, respectively. Suppose $\mathbb{P}$ and $\mathbb{Q}$ are equivalent measures satisfying Eq.~(\ref{emm_cts}). If $\lambda_{+}<2\tilde{\lambda}_{+}$ and $\lambda_{-}<2\tilde{\lambda}_{-}$, then we have
		\begin{align}
		\label{kld_cts}
			KL(\mathbb{P}||\mathbb{Q})=TC\Gamma(-a)\Big(\big((a-1)\lambda_+^{a}-a\tilde{\lambda}_{+}\lambda_+^{a-1}+\tilde{\lambda}_{+}^a\big)+\big((a-1)\lambda_-^{a}-a\tilde{\lambda}_{-}\lambda_-^{a-1}+\tilde{\lambda}_{-}^a\big)\Big).
		\end{align}
	\end{prp}
	
	For details of the proof, please check Theorem~2 and its proof in Kim and Lee \cite{kim2007relative}. Shortly, plugging Eq.~(\ref{log_radon_nikodym}) to Eq.~(\ref{kld_measure}) provides the proof.

	Note that some parts of the original theorem in Kim and Lee \cite{kim2007relative} are modified based on the parametrization and equivalent martingale measure condition given in this paper.
	
	As mentioned in Section~\ref{sec_intro}, the Kullback--Leibler divergence for GTS processes has not been found.
	
\subsection{RDTS processes}
	Another tempered stable process we cover in this paper is the RDTS process \cite{kim2010tempered,rachev2011financial}. 
	
	Similar to a GTS distribution, an RDTS distribution has seven parameters, $\boldsymbol{\xi}=(a_+, a_-, C_+, C_-,\lambda_+,\lambda_-,m)$: $m$ is the location parameter; $a_+$ and $a_-$ are the tail indexes; $C_+$ and $C_-$ are the scale parameters; $\lambda_+$ and $\lambda_-$ are the decay rates of the upper and lower tails, respectively. $C_+, C_-, \lambda_+,\lambda_-$ are positive, $a_+, a_-\in(0,2)\setminus \{1\}$, and $m\in\mathbb{R}$.
	
	An RDTS distribution has the following tempering function:
	 \begin{align}
	 \label{temp_fns_rdts}
	 	t_{\textrm{RDTS}}(x;\boldsymbol{\xi})={\rm e}^{-\frac{\lambda_{+}}{2}x^2} 1_{x>0}(x)+{\rm e}^{-\frac{\lambda_{-}}{2} |x|^2} 1_{x<0}(x),
	 \end{align}
	and the L\'evy measure of the RDTS distribution is obtained from Eq.~(\ref{levy_measure_ts}) \cite{rachev2011financial} as
	\begin{align}
	\label{levy_measure_rdts}
		\nu(dx;\boldsymbol{\xi})=\Big(\frac{C_+{\rm e}^{-\frac{\lambda_+}{2} x^2}}{x^{a_++1}} 1_{x>0}(x)+\frac{C_-{\rm e}^{-\frac{\lambda_-}{2} |x|^2}}{|x|^{a_-+1}} 1_{x<0}(x)\Big)dx.
	\end{align}
	
	From Eq.~(\ref{charac_ftn_ts}), the characteristic function of an RDTS distribution \cite{rachev2011financial} is represented with
	\begin{align}
	\label{charac_ftn_rdts}
		\phi(u;\boldsymbol{\xi})=&\exp\Big(ium+C_+G(iu;a_+,\lambda_+)+C_-G(-iu;a_-,\lambda_-)\Big),
	\end{align}
	where
	\begin{align}
	\begin{split}
		G(x;a,\lambda)=&2^{-1-\frac{a}{2}}\lambda^a\Gamma(-\frac{a}{2})\big(M(-\frac{a}{2},\frac{1}{2};\frac{x^2}{2\lambda^2})-1\big)\\
		&+2^{-\frac{1}{2}-\frac{a}{2}}x\lambda^{a-1}\Gamma(\frac{1-a}{2})\big(M(\frac{1-a}{2},\frac{3}{2};\frac{x^2}{2\lambda^2})-1\big),
	\end{split}
	\end{align}
	such that $M$ is the confluent hypergeometric function, and $\Gamma$ is the Gamma function. However, a closed-form expression for the probability density function of RDTS distributions does not exist.
		
	Similar to GTS processes and CTS processes, an RDTS process with parameters $\boldsymbol{\xi}=(a_+, a_-, C_+, C_-, \lambda_+, \lambda_-, m)$ is a L\'evy process with the L\'evy measure in Eq.~(\ref{levy_measure_rdts}) and the drift term determined by $m$.
	
	The equivalent martingale measure condition for RDTS processes \cite{rachev2011financial} is as follows:
	\begin{align}
	\label{emm_rdts}
	\left\{ 
		\begin{array}{ll}
			C_+=\tilde{C}_+\\
			C_-=\tilde{C}_-\\
			a_+=\tilde{a}_+\\
			a_-=\tilde{a}_-\\
			m=\tilde{m}+2^{-\frac{1+a_+}{2}}C_+\Gamma(\frac{1-a_+}{2})(\lambda_+^{a_+-1}-\tilde{\lambda}_{+}^{a_+-1})-2^{-\frac{1+a_-}{2}}C_-\Gamma(\frac{1-a_-}{2})(\lambda_-^{a_--1}-\tilde{\lambda}_{-}^{a_--1})
		\end{array}
		\right..
	\end{align}
	
	For equivalent martingale measures, the Radon--Nikodym derivative between the L\'evy measures of RDTS processes is obtained from Eq.~(\ref{levy_measure_rdts}) as
	\begin{align}
	\label{radon_nikodym_rdts}
		\frac{d\nu^\mathbb{P}}{d\nu^\mathbb{Q}}={\rm e}^{-\frac{(\lambda_+-\tilde{\lambda}_+)}{2}x^2} 1_{x>0}(x)+{\rm e}^{-\frac{(\lambda_{-}-\tilde{\lambda}_{-})}{2}|x|^2} 1_{x<0}(x).
	\end{align}
	
	The logarithmic Radon--Nikodym derivative for two equivalent martingale RDTS processes is expressed with
	\begin{align}
	\label{log_radon_nikodym_rdts}
		\psi(x;\boldsymbol{\xi})=-\frac{(\lambda_+-\tilde{\lambda}_+)}{2}x^2 1_{x>0}(x)-\frac{(\lambda_{-}-\tilde{\lambda}_{-})}{2}|x|^2 1_{x<0}(x).
	\end{align}
	
\section{Information geometry for tempered stable processes}
\label{sec_ig}
	In this section, we present the main results of this paper. By deriving the $\alpha$-divergence between two tempered stable processes, we explore the information geometry of tempered stable processes. The metric tensor and the $\alpha$-connection of the geometry are calculated for various tempered stable processes.
	 
	Our starting point is the $\alpha$-divergence between two probability density functions. The $\alpha$-divergence between two probability density functions $p$ and $q$ \cite{amari2000methods} is defined as 
	 \begin{align}
	 \label{alpha_div}
	D^{(\alpha )}(p||q)=\left\{ 
	\begin{array}{ll}
	\frac{4}{1-\alpha^{2}}\int (\frac{1-\alpha}{2}p(x;\boldsymbol{\xi})+\frac{1+\alpha}{2}q(x;\tilde{\boldsymbol{\xi}})-p(x;\boldsymbol{\xi})^{\frac{1-\alpha}{2}}q(x;\tilde{\boldsymbol{\xi}})^{\frac{1+\alpha}{2}})dx & (\alpha \neq \pm 1)\\ 
	\int (p(x;\boldsymbol{\xi}) \log{\frac{p(x;\boldsymbol{\xi})}{q(x;\tilde{\boldsymbol{\xi}})}}-p(x;\boldsymbol{\xi})+q(x;\tilde{\boldsymbol{\xi}}))dx & (\alpha =-1)\\
	\int (q(x;\tilde{\boldsymbol{\xi}}) \log{\frac{q(x;\tilde{\boldsymbol{\xi}})}{p(x;\boldsymbol{\xi})}}-q(x;\tilde{\boldsymbol{\xi}})+p(x;\boldsymbol{\xi}))dx & (\alpha =1)
	\end{array}
	\right.,
	\end{align}
	where $\boldsymbol{\xi}$ and $\tilde{\boldsymbol{\xi}}$ are parameters of $p$ and $q$, respectively. It is straightforward to check that the $\alpha$-divergence is not symmetric under the exchange between $p$ and $q$ except for $\alpha=0$.
	
	It is noteworthy that several $\alpha$ values are related to well-known divergences and distances. $D^{(-1)}(p||q)$ is the Kullback--Leibler divergence, also known as relative entropy. The Hellinger distance corresponds to $\alpha=0$.
	
	The $\alpha$-divergence has an interesting characteristic called $\alpha$-duality \cite{amari2000methods}: changing the sign of $\alpha$ corresponds to the exchange between $p$ and $q$ such that
	\begin{align}
	\label{a_duality}
		D^{(\alpha )}(p||q)=D^{(-\alpha )}(q||p).
	\end{align}
	It is straightforward to check that the Hellinger distance, which is the $0$-divergence, is self-dual and symmetric in the exchange of $p$ and $q$. $D^{(1)}(p||q)$ is the dual to the Kullback--Leibler divergence.
			
	As given in Eq.~(\ref{kld_measure}), the Kullback--Leibler divergence for tempered stable processes is given \cite{cont2004nonparametric,kim2007relative} as
	\begin{align}
	\label{kld_measure_ts}
		KL(\mathbb{P}||\mathbb{Q})=T \int \Big(\Big(\frac{d\nu^\mathbb{P}}{d\nu^\mathbb{Q}}\Big) \log{\Big(\frac{d\nu^\mathbb{P}}{d\nu^\mathbb{Q}}\Big)}- \Big(\frac{d\nu^\mathbb{P}}{d\nu^\mathbb{Q}}\Big)+1\Big) \nu^\mathbb{Q}(dx).
	\end{align}
	
	For exploring the information geometry, let us extend Eq.~(\ref{kld_measure_ts}) based on the Kullback--Leibler divergence to $\alpha$-divergence by using $f$-divergence. 

	For equivalent pure-jump L\'evy processes, the Radon--Nikodym derivative between the path measures $\mathbb{P}$ and $\mathbb{Q}$ can be represented in terms of the logarithmic Radon--Nikodym derivative between their L\'evy measures in Eq.~(\ref{log_radon_nikodym}) \cite{sato1999levy,sato2000density}. In particular, the logarithmic Radon--Nikodym derivative is given by
	\begin{align}
	\label{rn_path}
	\begin{split}
	\log\Big(\frac{d\mathbb{P}}{d\mathbb{Q}}\Big)	&=\int_0^T\int \psi(x)\left(N(dt,dx)-dt\,\nu^\mathbb{Q}(dx)\right)+T\int\left(\psi(x)-{\rm e}^{\psi(x)}+1\right)\nu^\mathbb{Q}(dx)\\
	&=\int_0^T\int\psi(x)N(dt,dx)-T\int\left({\rm e}^{\psi(x)}-1\right)\nu^\mathbb{Q}(dx),
	\end{split}
	\end{align}
	where $N(dt,dx)$ is the jump measure under $\mathbb{Q}$, defined by $N((0,t]\times A)=\sum_{0<s\leq t}1_A(\Delta X_s)$ for a Borel set $A$ bounded away from zero. 
	
	For $0<r<1$, the $r$-th moment of the Radon--Nikodym derivative $d\mathbb{P}/d\mathbb{Q}$ under $\mathbb{Q}$ is known as the Hellinger integral \cite{sato1999levy,sato2000density}:
	\begin{align}
	\label{hellinger_levy}
	H_r(\mathbb{P},\mathbb{Q})=\mathbb{E}_{\mathbb{Q}}\left[\left(\frac{d\mathbb{P}}{d\mathbb{Q}}\right)^r\right]=\exp\left[-T\int\left(r{\rm e}^{\psi(x)}+1-r-{\rm e}^{r\psi(x)}\right)\nu^\mathbb{Q}(dx)\right].
	\end{align}
	The same expression for $r>1$ follows from the positive-order R\'enyi-divergence formula for Poisson random measures \cite{leskela2024information}.

	\begin{thm}
	\label{thm_div_ts}
	Let $(X_t, \mathbb{P})_{t\in[0,T]}$ and $(X_t, \mathbb{Q})_{t\in[0,T]}$ be tempered stable processes with L\'evy triplets of $(0,\nu^\mathbb{P},\gamma^{\mathbb{P}})$ and $(0,\nu^\mathbb{Q},\gamma^{\mathbb{Q}})$, respectively. Suppose $\mathbb{P}$ and $\mathbb{Q}$ are equivalent measures. Then we have the $\alpha$-divergence between two tempered stable processes as
	 \begin{align}
	 \label{a_div_measure_ts}
	D^{(\alpha )}(\mathbb{P}||\mathbb{Q})=\left\{ 
	\begin{array}{ll}
	\frac{4}{1-\alpha ^{2}}\left\{1-\exp\left[-T\int \left(\frac{1-\alpha}{2}\left(\frac{d\nu^\mathbb{P}}{d\nu^\mathbb{Q}}\right) +\frac{1+\alpha}{2}-\left(\frac{d\nu^\mathbb{P}}{d\nu^\mathbb{Q}}\right)^{\frac{1-\alpha}{2}}\right) \nu^\mathbb{Q}(dx)\right]\right\} & (\alpha \neq \pm 1)\\ 
	T\int \left(\left(\frac{d\nu^\mathbb{P}}{d\nu^\mathbb{Q}}\right) \log{\left(\frac{d\nu^\mathbb{P}}{d\nu^\mathbb{Q}}\right)}- \left(\frac{d\nu^\mathbb{P}}{d\nu^\mathbb{Q}}\right)+1\right) \nu^\mathbb{Q}(dx) & (\alpha =-1)\\
	T\int \left(\left(\frac{d\nu^\mathbb{Q}}{d\nu^\mathbb{P}}\right) \log{\left(\frac{d\nu^\mathbb{Q}}{d\nu^\mathbb{P}}\right)}- \left(\frac{d\nu^\mathbb{Q}}{d\nu^\mathbb{P}}\right)+1\right) \nu^\mathbb{P}(dx) & (\alpha =1)
	\end{array}
	\right..
	\end{align}
	\end{thm}
	\begin{proof}
	There is the well-known $f$-divergence which produces $\alpha$-divergence \cite{amari2000methods}:
	\begin{align}
		D^{(\alpha )}(\mathbb{P}||\mathbb{Q})=D_{f^{(\alpha)}}(\mathbb{P}||\mathbb{Q})=\mathbb{E}_{\mathbb{Q}}\left[f^{(\alpha)}\left(\frac{d\mathbb{P}}{d\mathbb{Q}}\right)\right],
	\end{align}
	where $f^{(\alpha)}(t)$ is given by
	\begin{align}
		f^{(\alpha)}(t)=\left\{ 
	\begin{array}{ll}
	\frac{4}{1-\alpha ^{2}}\Big(\frac{1-\alpha}{2} t +\frac{1+\alpha}{2}-t^{\frac{1-\alpha}{2}}\Big) & (\alpha \neq \pm 1)\\ 
	t \log{t} - t +1 & (\alpha =-1)\\
	- \log{t} +t - 1 & (\alpha =1)
	\end{array}
	\right..
	\end{align}
	Since $\mathbb{E}_{\mathbb{Q}}[d\mathbb{P}/d\mathbb{Q}]=1$, substituting the Radon--Nikodym derivative into the $\alpha$-divergence for $\alpha\neq\pm1$ gives
	\begin{align}
	D^{(\alpha)}(\mathbb{P}||\mathbb{Q})
	=\frac{4}{1-\alpha^2}\left\{1-\mathbb{E}_{\mathbb{Q}}\left[\left(\frac{d\mathbb{P}}{d\mathbb{Q}}\right)^{\frac{1-\alpha}{2}}\right]\right\}.
	\end{align}
	Let us set $r=\frac{1-\alpha}{2}$. For $|\alpha|<1$, i.e., $0<r<1$, the $\alpha$-divergence from Eq.~(\ref{hellinger_levy}) is given by
	\begin{align}
	\label{a_div_neq_pm1}
	\begin{split}
	&D^{(\alpha)}(\mathbb{P}||\mathbb{Q})	=\frac{4}{1-\alpha^2}\left\{1-H_{\frac{1-\alpha}{2}}(\mathbb{P},\mathbb{Q})\right\}\\
	&=\frac{4}{1-\alpha^2}\left\{1-\exp\left[-T\int\left(\frac{1-\alpha}{2}\frac{d\nu^\mathbb{P}}{d\nu^\mathbb{Q}}+\frac{1+\alpha}{2}-\left(\frac{d\nu^\mathbb{P}}{d\nu^\mathbb{Q}}\right)^{\frac{1-\alpha}{2}}\right)\nu^\mathbb{Q}(dx)\right]\right\}.
	\end{split}
	\end{align}
	For $\alpha<-1$, i.e., $r>1$, we obtain the same expression as Eq.~(\ref{a_div_neq_pm1}) by using Eq.~(\ref{hellinger_levy}) based on \cite{leskela2024information}. For $\alpha>1$, i.e., $r<0$, the result follows from the $\alpha$-duality of Eq.~(\ref{a_duality}) by exchanging $\mathbb{P}$ and $\mathbb{Q}$. From these, we derive Eq.~(\ref{a_div_measure_ts}) for $\alpha\neq\pm1$. 
		
	Applying L'Hopital's rule to the case $\alpha\neq\pm1$ in the limits $\alpha\to-1$ and $\alpha\to1$ yields the corresponding Kullback--Leibler divergences for $\alpha=-1$ and $\alpha=1$, respectively.
	\end{proof}

	The $\alpha$-divergence between two tempered stable processes can be represented with $\lambda(x;\boldsymbol{\xi})$ in Eq.~(\ref{levy_measure_ts}):
	\begin{align}
	 \label{a_div_measure_ts_lambda}
	D^{(\alpha )}(\mathbb{P}||\mathbb{Q})=\left\{ 
	\begin{array}{ll}
	\frac{4}{1-\alpha ^{2}}\left\{1-\exp\left[-T\int \left(\frac{1-\alpha}{2}\frac{\lambda^\mathbb{P}}{\lambda^\mathbb{Q}}+\frac{1+\alpha}{2}-\left(\frac{\lambda^\mathbb{P}}{\lambda^\mathbb{Q}}\right)^{\frac{1-\alpha}{2}}\right)\lambda^\mathbb{Q}dx\right]\right\} & (\alpha \neq \pm 1)\\ 
	T\int \Big(\Big(\frac{\lambda^\mathbb{P}}{\lambda^\mathbb{Q}}\Big) \log{\Big(\frac{\lambda^\mathbb{P}}{\lambda^\mathbb{Q}}\Big)}- \Big(\frac{\lambda^\mathbb{P}}{\lambda^\mathbb{Q}}\Big)+1\Big) \lambda^\mathbb{Q}dx & (\alpha =-1)\\
	T\int \Big(\Big(\frac{\lambda^\mathbb{Q}}{\lambda^\mathbb{P}}\Big) \log{\Big(\frac{\lambda^\mathbb{Q}}{\lambda^\mathbb{P}}\Big)}- \Big(\frac{\lambda^\mathbb{Q}}{\lambda^\mathbb{P}}\Big)+1\Big) \lambda^\mathbb{P}dx & (\alpha =1)
	\end{array}
	\right..
	\end{align}

	It is noteworthy that Eq.~(\ref{a_div_measure_ts_lambda}) can be represented only with the tempering function $t(x;\boldsymbol{\xi})$ in Eq.~(\ref{levy_measure_ts}) if the equivalent martingale measure condition on a tempered stable distribution includes the parameters in the L\'evy measure of the stable distribution, Eq.~(\ref{levy_measure_stable}), i.e., $C$ and $a$. In this case, since the Radon--Nikodym derivative can be expressed with the tempering function $t$:
	\begin{align}
		\frac{d\nu^\mathbb{P}}{d\nu^\mathbb{Q}}=\frac{\lambda^\mathbb{P}}{\lambda^\mathbb{Q}}=\frac{t^\mathbb{P}}{t^\mathbb{Q}},
	\end{align}
	the $\alpha$-divergence is also given in terms of the tempering function $t$.
			
	Additionally, the $\alpha$-divergence between two tempered stable processes, Eq.~(\ref{a_div_measure_ts}), can be expressed in terms of the logarithmic Radon--Nikodym derivative, Eq.~(\ref{log_radon_nikodym}):
	\begin{align}
	\label{a_div_measure_ts_psi}
	D^{(\alpha )}(\mathbb{P}||\mathbb{Q})=\left\{ 
	\begin{array}{ll}
	\frac{4}{1-\alpha ^{2}}\left\{1-\exp\left[-T\int \left(\frac{1-\alpha}{2}{\rm e}^{\psi(x)}+\frac{1+\alpha}{2}-{\rm e}^{\frac{1-\alpha}{2}\psi(x)}\right)\nu^\mathbb{Q}(dx)\right]\right\} & (\alpha \neq \pm 1)\\ 
	T\int (\psi(x){\rm e}^{\psi(x)}-{\rm e}^{\psi(x)}+1)\nu^\mathbb{Q}(dx) & (\alpha =-1)\\
	T\int (-\psi(x){\rm e}^{-\psi(x)}-{\rm e}^{-\psi(x)}+1)\nu^\mathbb{P}(dx) & (\alpha =1)
	\end{array}
	\right..
	\end{align}

	It is straightforward to derive the information geometry from the $\alpha$-divergence given in Theorem~\ref{thm_div_ts}. The following theorem gives the information geometry of tempered stable processes.
	\begin{thm}
	\label{thm_geo_ts}
	Let $(X_t, \mathbb{P})_{t\in[0,T]}$ be a tempered stable process with the L\'evy triplet of $(0,\nu,\gamma)$. The metric tensor and the $\alpha$-connection for the information geometry of tempered stable processes are given by
	\begin{align}
		\label{metric_ts}
		g_{ij}=&T\int \partial_i \log{\Big(\frac{d\nu}{dx}\Big)} \partial_j \log{\Big(\frac{d\nu}{dx}\Big)} \nu(dx),\\
		\label{conn_ts}
		\Gamma^{(\alpha)}_{ij,k}=&T\int \Big(\partial_i \partial_j \log{\Big(\frac{d\nu}{dx}\Big)}+\frac{1-\alpha}{2}\partial_i \log{\Big(\frac{d\nu}{dx}\Big)} \partial_j \log{\Big(\frac{d\nu}{dx}\Big)}\Big)\partial_k \log{\Big(\frac{d\nu}{dx}\Big)} \nu(dx),
	\end{align}
	where $i,j,$ and $k$ run for the coordinate system $\boldsymbol{\xi}$.
	\end{thm}
	\begin{proof}
	For a given divergence $D$, the metric tensor and the connection of the geometry \cite{amari2000methods} are derived as
	\begin{align}
		\label{ig_metric}
		g_{ij}&=-D(\partial_i, \tilde{\partial}_j)|_{\boldsymbol{\xi}=\tilde{\boldsymbol{\xi}}},\\
		\label{ig_connection}
		\Gamma_{ij,k}&=-D(\partial_i \partial_j,\tilde{\partial}_k)|_{\boldsymbol{\xi}=\tilde{\boldsymbol{\xi}}},
	\end{align}
	where $i,j,$ and $k$ run for the coordinate system $\boldsymbol{\xi}$.
	
	By plugging the $\alpha$-divergence of Eq.~(\ref{a_div_measure_ts}) to Eq.~(\ref{ig_metric}) and Eq.~(\ref{ig_connection}), we easily obtain the following:
	\begin{align}
		g_{ij}=&T\int \partial_i \log{\Big(\frac{d\nu}{dx}\Big)} \partial_j \log{\Big(\frac{d\nu}{dx}\Big)} \nu(dx),\nonumber\\
		\Gamma^{(\alpha)}_{ij,k}=&T\int \Big(\partial_i \partial_j \log{\Big(\frac{d\nu}{dx}\Big)}+\frac{1-\alpha}{2}\partial_i \log{\Big(\frac{d\nu}{dx}\Big)} \partial_j \log{\Big(\frac{d\nu}{dx}\Big)}\Big)\partial_k \log{\Big(\frac{d\nu}{dx}\Big)} \nu(dx),\nonumber
	\end{align}
	where $i,j,$ and $k$ run for the coordinate system $\boldsymbol{\xi}$.
	\end{proof}
	
	It is noteworthy that the metric tensor of information geometry is the Fisher information matrix, i.e., Eq.~(\ref{metric_ts}) is the Fisher information matrix of tempered stable processes. Additionally, Eq.~(\ref{metric_ts}) and Eq.~(\ref{conn_ts}) have the same form as the corresponding Fisher information matrix and the $\alpha$-connection obtained from the probability density function $p$.
	
	Similar to the $\alpha$-divergence, both the metric tensor and the $\alpha$-connection of information geometry are expressed in terms of $\lambda(x;\boldsymbol{\xi})$ in Eq.~(\ref{levy_measure_ts}):
	\begin{align}
	\label{metric_ts_lambda}
		g_{ij}=&T\int \partial_i \log{\lambda(x;\boldsymbol{\xi})} \partial_j \log{\lambda(x;\boldsymbol{\xi})} \lambda(x;\boldsymbol{\xi}) dx,\\
	\label{conn_ts_lambda}
		\Gamma^{(\alpha)}_{ij,k}=&T\int \Big(\partial_i \partial_j \log{\lambda(x;\boldsymbol{\xi})}+\frac{1-\alpha}{2}\partial_i \log{\lambda(x;\boldsymbol{\xi})} \partial_j \log{\lambda(x;\boldsymbol{\xi})}\Big)\Big(\partial_k \log{\lambda(x;\boldsymbol{\xi})}\Big) \lambda(x;\boldsymbol{\xi}) dx,
	\end{align}
	where $i,j,$ and $k$ run for the coordinate system $\boldsymbol{\xi}$. It is also noteworthy that Eq.~(\ref{metric_ts_lambda}) and Eq.~(\ref{conn_ts_lambda}) have the same form as those using the probability density function given in Amari and Nagaoka \cite{amari2000methods}.
	
	If the equivalent martingale measure condition on a tempered stable distribution includes the parameters in the L\'evy measure of the stable distribution in Eq.~(\ref{levy_measure_stable}), the metric tensor and the $\alpha$-connection are expressed with the L\'evy measure and the tempering function $t(x;\boldsymbol{\xi})$ in Eq.~(\ref{levy_measure_ts}):
	\begin{align}
	\label{metric_ts_t}
		g_{ij}=&T\int \partial_i \log{t(x;\boldsymbol{\xi})} \partial_j \log{t(x;\boldsymbol{\xi})} \nu(dx),\\
	\label{conn_ts_t}
		\Gamma^{(\alpha)}_{ij,k}=&T\int \Big(\partial_i \partial_j \log{t(x;\boldsymbol{\xi})}+\frac{1-\alpha}{2}\partial_i \log{t(x;\boldsymbol{\xi})} \partial_j \log{t(x;\boldsymbol{\xi})}\Big)\Big(\partial_k \log{t(x;\boldsymbol{\xi})}\Big) \nu(dx),
	\end{align}
	where $i,j,$ and $k$ run for the coordinate system $\boldsymbol{\xi}$.
	
\subsection{GTS processes}
	Theorem~\ref{thm_div_ts} can be applied to GTS processes. By plugging Eq.~(\ref{radon_nikodym_gts}) to Eq.~(\ref{a_div_measure_ts}) in Theorem~\ref{thm_div_ts}, the $\alpha$-divergence of GTS processes is obtained as given below.
	
	\begin{crl}
	\label{crl_div_gts}
		Let $(X_t, \mathbb{P})_{t\in[0,T]}$ and $(X_t, \mathbb{Q})_{t\in[0,T]}$ be GTS processes with parameters $\boldsymbol{\xi}=(a_+, a_-, C_+, C_-,\lambda_+,\lambda_-,m)$ and $\tilde{\boldsymbol{\xi}}=(\tilde{a}_+, \tilde{a}_-, \tilde{C}_+, \tilde{C}_-,\tilde{\lambda}_+,\tilde{\lambda}_-,\tilde{m})$, respectively. Suppose $\mathbb{P}$ and $\mathbb{Q}$ are equivalent measures satisfying Eq.~(\ref{emm_gts}). Then we have the $\alpha$-divergence between two GTS processes as
	\begin{align}
	\label{a_div_gts}
	D^{(\alpha )}(\mathbb{P}||\mathbb{Q})=\left\{ 
	\begin{array}{ll}
	\begin{aligned}
\frac{4}{1-\alpha^2}\Bigg\{1-\exp\Bigg[&-TC_+\Gamma(-a_+)\Big(\frac{1-\alpha}{2}\lambda_+^{a_+}+\frac{1+\alpha}{2}\tilde{\lambda}_+^{a_+}\\
&\qquad-\big(\frac{1-\alpha}{2}\lambda_++\frac{1+\alpha}{2}\tilde{\lambda}_+\big)^{a_+}\Big)\\
&-TC_-\Gamma(-a_-)\Big(\frac{1-\alpha}{2}\lambda_-^{a_-}+\frac{1+\alpha}{2}\tilde{\lambda}_-^{a_-}\\
&\qquad-\big(\frac{1-\alpha}{2}\lambda_-+\frac{1+\alpha}{2}\tilde{\lambda}_-\big)^{a_-}\Big)\Bigg]\Bigg\}
\end{aligned} & (\alpha \neq \pm 1)\\ 
	TC_+\Gamma(-a_+)\big((a_+-1)\lambda_+^{a_+}-a_+\tilde{\lambda}_+\lambda_+^{a_+-1}+\tilde{\lambda}_+^{a_+}\big)\\
	+TC_-\Gamma(-a_-)\big((a_--1)\lambda_-^{a_-}-a_-\tilde{\lambda}_-\lambda_-^{a_--1}+\tilde{\lambda}_-^{a_-}\big) & (\alpha =-1)\\
	TC_+\Gamma(-a_+)\big((a_+-1)\tilde{\lambda}_+^{a_+}-a_+\lambda_+\tilde{\lambda}_+^{a_+-1}+\lambda_+^{a_+}\big)\\
	+TC_-\Gamma(-a_-)\big((a_--1)\tilde{\lambda}_-^{a_-}-a_-\lambda_-\tilde{\lambda}_-^{a_--1}+\lambda_-^{a_-}\big) & (\alpha =1)
	\end{array}
	\right..
	\end{align}
	\end{crl}
	\begin{proof}
	We start with $\alpha\neq\pm1$. From Eq.~(\ref{hellinger_levy}) and Eq.~(\ref{a_div_measure_ts}) in Theorem~\ref{thm_div_ts}, the positive $x$ term of $-\log H_{\frac{1-\alpha}{2}}(\mathbb{P},\mathbb{Q})$ is given by
	\begin{align}
		&-\log H_{\frac{1-\alpha}{2}}(\mathbb{P},\mathbb{Q})|_{+}=T\int_0^\infty \Big(\frac{1-\alpha}{2}{\rm e}^{-(\lambda_+-\tilde{\lambda}_+)x} +\frac{1+\alpha}{2}- \Big({\rm e}^{-(\lambda_+-\tilde{\lambda}_+)x}\Big)^{\frac{1-\alpha}{2}}\Big)\frac{C_+{\rm e}^{-\tilde{\lambda}_+ x}}{x^{a_++1}} dx\nonumber\\
	&=TC_+\int_0^\infty \Big(\frac{1-\alpha}{2}\frac{{\rm e}^{-\lambda_+ x}}{x^{a_++1}} +\frac{1+\alpha}{2}\frac{{\rm e}^{-\tilde{\lambda}_+ x}}{x^{a_++1}} -\frac{{\rm e}^{-(\frac{1-\alpha}{2}\lambda_++\frac{1+\alpha}{2}\tilde{\lambda}_+) x}}{x^{a_++1}} \Big)dx\nonumber\\
	&=TC_+\int_0^\infty \Big(\frac{1-\alpha}{2}\frac{\lambda_+^{a_+}{\rm e}^{-t}}{t^{a_++1}} +\frac{1+\alpha}{2}\frac{\tilde{\lambda}_+^{a_+}{\rm e}^{-t}}{t^{a_++1}} -\frac{(\frac{1-\alpha}{2}\lambda_++\frac{1+\alpha}{2}\tilde{\lambda}_+)^{a_+}{\rm e}^{-t}}{t^{a_++1}} \Big)dt\nonumber\\
	&=TC_+\Gamma(-a_+)\Big(\frac{1-\alpha}{2}\lambda_+^{a_+}+\frac{1+\alpha}{2}\tilde{\lambda}_+^{a_+}-\big(\frac{1-\alpha}{2}\lambda_++\frac{1+\alpha}{2}\tilde{\lambda}_+\big)^{a_+}\Big).\nonumber
	\end{align}
	After a similar calculation is done for the negative $x$ term, we obtain
	\begin{align}
	-\log H_{\frac{1-\alpha}{2}}(\mathbb{P},\mathbb{Q})=&TC_+\Gamma(-a_+)\Big(\frac{1-\alpha}{2}\lambda_+^{a_+}+\frac{1+\alpha}{2}\tilde{\lambda}_+^{a_+}-\big(\frac{1-\alpha}{2}\lambda_++\frac{1+\alpha}{2}\tilde{\lambda}_+\big)^{a_+}\Big)\nonumber\\
	&+TC_-\Gamma(-a_-)\Big(\frac{1-\alpha}{2}\lambda_-^{a_-}+\frac{1+\alpha}{2}\tilde{\lambda}_-^{a_-}-\big(\frac{1-\alpha}{2}\lambda_-+\frac{1+\alpha}{2}\tilde{\lambda}_-\big)^{a_-}\Big).\nonumber
	\end{align}
	Substituting this expression into Eq.~(\ref{a_div_measure_ts}) gives the $\alpha$-divergence.

	For $\alpha=-1$, there are two ways to prove the result. The first proof is based on the direct computation by plugging the Radon--Nikodym derivative of Eq.~(\ref{radon_nikodym_gts}) to the $\alpha=-1$ case in Eq.~(\ref{a_div_measure_ts}). The second proof applies L'Hopital's rule to $\alpha \neq \pm1$ in the limit of $\alpha\to-1$. Both methods give the same result.
	
	We show the first proof here. Similar to $\alpha \neq \pm1$, let us start with the positive $x$ term of the Kullback--Leibler divergence:
	\begin{align}
	D^{(-1)}(\mathbb{P}||\mathbb{Q})|_{+}&=TC_+\int_0^\infty \Big(\frac{{\rm e}^{-\lambda_+x}}{x^{a_++1}}\big(-(\lambda_+-\tilde{\lambda}_+)x-1\big) + \frac{{\rm e}^{-\tilde{\lambda}_+x}}{x^{a_++1}}\Big)dx\nonumber\\
	&=TC_+\Gamma(-a_+)\big(\lambda_+^{a_+-1}a_+(\lambda_+-\tilde{\lambda}_+)-\lambda_+^{a_+}+\tilde{\lambda}_+^{a_+}\big)\nonumber\\
	&=TC_+\Gamma(-a_+)\big((a_+-1)\lambda_+^{a_+}-a_+\tilde{\lambda}_+\lambda_+^{a_+-1}+\tilde{\lambda}_+^{a_+}\big).\nonumber
	\end{align}
	
	After repeating the same step for negative $x$, the Kullback--Leibler divergence between two GTS processes is given by
	\begin{align}
	D^{(-1)}(\mathbb{P}||\mathbb{Q})=&TC_+\Gamma(-a_+)\big((a_+-1)\lambda_+^{a_+}-a_+\tilde{\lambda}_+\lambda_+^{a_+-1}+\tilde{\lambda}_+^{a_+}\big)\nonumber\\
	&+TC_-\Gamma(-a_-)\big((a_--1)\lambda_-^{a_-}-a_-\tilde{\lambda}_-\lambda_-^{a_--1}+\tilde{\lambda}_-^{a_-}\big).\nonumber
	\end{align}
	It is trivial that we obtain the same result given above when applying L'Hopital's rule to the $\alpha \neq \pm1$ case.
	
	It is possible to use the same ways mentioned above in order to derive the $1$-divergence. Alternatively, we also can use $\alpha$-duality with $\alpha=-1$. By following a similar calculation, the $\alpha$-divergence is obtained as
	\begin{align}
	D^{(1)}(\mathbb{P}||\mathbb{Q})=&D^{(-1)}(\mathbb{Q}||\mathbb{P})\nonumber\\
	=&TC_+\Gamma(-a_+)\big((a_+-1)\tilde{\lambda}_+^{a_+}-a_+\lambda_+\tilde{\lambda}_+^{a_+-1}+\lambda_+^{a_+}\big)\nonumber\\
	&+TC_-\Gamma(-a_-)\big((a_--1)\tilde{\lambda}_-^{a_-}-a_-\lambda_-\tilde{\lambda}_-^{a_--1}+\lambda_-^{a_-}\big).\nonumber
	\end{align}
	It is also straightforward to obtain the same result by using the direct calculation or L'Hopital's rule.
	\end{proof}

	In another way, we can derive the $\alpha$-divergence between GTS processes from plugging Eq.~(\ref{log_radon_nikodym_gts}) to Eq.~(\ref{a_div_measure_ts_psi}). Additionally, the $\alpha$-divergence can be obtained from Eq.~(\ref{a_div_measure_ts_lambda}) or a similar expression using the tempering function of GTS processes.
	
	We next derive the information-geometric quantities of GTS processes from Theorem~\ref{thm_geo_ts}.
	
	It is noteworthy that although GTS processes have seven parameters, their geometry is reduced to a two-dimensional manifold due to the equivalent martingale measure condition of Eq.~(\ref{emm_gts}). The coordinate system for the GTS geometry is simply $\boldsymbol{\xi}=(\lambda_+,\lambda_-)$.
	
	From Eq.~(\ref{levy_measure_gts}) and Eq.~(\ref{metric_ts}) in Theorem~\ref{thm_geo_ts}, the Fisher information matrix of GTS processes is found as
	\begin{align}
	\label{metric_gts}
		g_{ij}=\Bigg(
		\begin{array}{cc}
			\frac{TC_+\Gamma(2-a_+)}{\lambda_+^{2-a_+}} &0 \\ 
			0 & \frac{TC_-\Gamma(2-a_-)}{\lambda_-^{2-a_-}}
		\end{array}\Bigg),
	\end{align}
	where $i, j$ run for $\lambda_+$ and $\lambda_-$. We can obtain the same result from Eq.~(\ref{metric_ts_lambda}) with the L\'evy measure or Eq.~(\ref{metric_ts_t}) with the GTS tempering function, Eq.~(\ref{temp_fns_gts}).
	
	It is easy to verify that the metric tensor is diagonal. This diagonality makes sense due to the fact that the L\'evy measure of GTS processes, Eq.~(\ref{levy_measure_gts}), is split by the $\lambda_+$ part and the $\lambda_-$ part using the indicator function. Additionally, $2-a_{\pm}$ are positive based on the parametrization of the GTS distribution.
	
	Since the Levi-Civita connection is defined as
	\begin{align}
	\label{lc_conn}
		\Gamma^{LC}_{ij,k}=\frac{1}{2}(\partial_i g_{jk} +\partial_j g_{ik}-\partial_k g_{ij}),
	\end{align}
	the non-trivial components of the Levi-Civita connection in GTS geometry are obtained as
	\begin{align}
	\label{lc_conn_gts_p}
		\Gamma^{LC}_{\lambda_+\lambda_+,\lambda_+}&=-\frac{1}{2}\frac{TC_+\Gamma(3-a_+)}{\lambda_+^{3-a_+}},\\
	\label{lc_conn_gts_m}
		\Gamma^{LC}_{\lambda_-\lambda_-,\lambda_-}&=-\frac{1}{2}\frac{TC_-\Gamma(3-a_-)}{\lambda_-^{3-a_-}}.
	\end{align}

	From Eq.~(\ref{levy_measure_gts}) and Eq.~(\ref{conn_ts}) in Theorem~\ref{thm_geo_ts}, the non-trivial $\alpha$-connection components of the GTS manifolds are given by
	\begin{align}
	\label{a_conn_gts_p}
	\Gamma^{(\alpha)}_{\lambda_+\lambda_+,\lambda_+}&=
	-\frac{1-\alpha}{2}\frac{TC_+\Gamma(3-a_+)}{ \lambda_+^{3-a_+}},\\
	\label{a_conn_gts_m}
	\Gamma^{(\alpha)}_{\lambda_-\lambda_-,\lambda_-}&=
	-\frac{1-\alpha}{2}\frac{TC_-\Gamma(3-a_-)}{ \lambda_-^{3-a_-}}.
	\end{align}
	All other components vanish. By plugging $\alpha=0$ to Eq.~(\ref{a_conn_gts_p}) and Eq.~(\ref{a_conn_gts_m}), we can obtain the same Levi-Civita connection, Eq.~(\ref{lc_conn_gts_p}) and Eq.~(\ref{lc_conn_gts_m}), respectively. This is expected because the $0$-connection is the Levi-Civita connection. The same result can be calculated from Eq.~(\ref{conn_ts_lambda}) with the L\'evy measure or Eq.~(\ref{conn_ts_t}) with the GTS tempering function of Eq.~(\ref{temp_fns_gts}).

	We further investigate the dual flatness of the GTS geometry. For $\alpha=1$, Eq.~(\ref{a_conn_gts_p}) and Eq.~(\ref{a_conn_gts_m}) vanish. The coordinates $(\lambda_+,\lambda_-)$ are affine coordinates for the $1$-connection, and the GTS geometry is e-flat.

	A potential function for the Fisher metric is obtained as
	\begin{align}
	\label{potential_gts}
	\Phi_{\textrm{GTS}}(\boldsymbol{\xi})
	&=TC_+\Gamma(-a_+)\lambda_+^{a_+}+TC_-\Gamma(-a_-)\lambda_-^{a_-}.
	\end{align}
	Using $\Gamma(2-a)=a(a-1)\Gamma(-a)$, the Hessian of Eq.~(\ref{potential_gts}) satisfies
	\begin{align}
	\partial_{\lambda_\pm}^2\Phi_{\textrm{GTS}}
	&=TC_\pm\Gamma(2-a_\pm)\lambda_\pm^{a_\pm-2}
	=g_{\lambda_\pm\lambda_\pm},
	\end{align}
	and the mixed derivatives vanish. These are exactly the components of the Fisher information matrix in Eq.~(\ref{metric_gts}).

	From the potential function, the dual coordinates are obtained as
	\begin{align}
	\label{dual_coord_gts}
	\eta_\pm&=\partial_{\lambda_\pm}\Phi_{\textrm{GTS}}
	=TC_\pm a_\pm\Gamma(-a_\pm)\lambda_\pm^{a_\pm-1}.
	\end{align}
	After raising the last index in Eq.~(\ref{a_conn_gts_p}) and Eq.~(\ref{a_conn_gts_m}), we obtain
	\begin{align}
		\Gamma^{(\alpha)\,\lambda_\pm}_{\lambda_\pm\lambda_\pm}=-\frac{1-\alpha}{2}\frac{2-a_\pm}{\lambda_\pm}.
	\end{align}
	For $\alpha=-1$, the connections are given by
	\begin{align}
		\Gamma^{(-1)\,\lambda_\pm}_{\lambda_\pm\lambda_\pm}=-\frac{2-a_\pm}{\lambda_\pm}.
	\end{align}
	Since $\partial_{\lambda_\pm}\eta_\pm=TC_\pm\Gamma(2-a_\pm)\lambda_\pm^{a_\pm-2}$, the dual coordinates in Eq.~(\ref{dual_coord_gts}) satisfy
	\begin{align}
		(\nabla^{(-1)}d\eta_\pm)_{\lambda_\pm\lambda_\pm}=\partial_{\lambda_\pm}^2\eta_\pm-\Gamma^{(-1)\,\lambda_\pm}_{\lambda_\pm\lambda_\pm}\partial_{\lambda_\pm}\eta_\pm=0.
	\end{align}
	The mixed components also vanish. We have $\nabla^{(-1)}d\eta_\pm=0$. In the coordinate system $(\eta_+,\eta_-)$, this is expressed as
	\begin{align}
		(\nabla^{(-1)}d\eta_k)_{\eta_i\eta_j}
		=-\widetilde{\Gamma}^{(-1)\,\eta_k}_{\eta_i\eta_j}=0,
	\end{align}
	where $i,j,k$ run for $+$ and $-$. The connection coefficients vanish in the dual coordinates. The coordinates $(\eta_+,\eta_-)$ are affine coordinates for the $(-1)$-connection. The $1$- and $(-1)$-connections are dual with respect to the Fisher metric, and both are flat. The GTS statistical manifold is dually flat.

	From the potential function, the Bregman divergence is defined as
	\begin{align}
		&B_{\Phi}(\boldsymbol{\xi},\tilde{\boldsymbol{\xi}})
		=\Phi(\boldsymbol{\xi})-\Phi(\tilde{\boldsymbol{\xi}})-\nabla\Phi(\tilde{\boldsymbol{\xi}})\cdot(\boldsymbol{\xi}-\tilde{\boldsymbol{\xi}}).
	\end{align}
	Using Eq.~(\ref{potential_gts}), the Bregman divergences for GTS processes are given by
	\begin{align}
		B_{\Phi}(\boldsymbol{\xi},\tilde{\boldsymbol{\xi}})
		&=TC_+\Gamma(-a_+)\big((a_+-1)\tilde{\lambda}_+^{a_+}-a_+\lambda_+\tilde{\lambda}_+^{a_+-1}+\lambda_+^{a_+}\big)\nonumber\\
		&\quad+TC_-\Gamma(-a_-)\big((a_--1)\tilde{\lambda}_-^{a_-}-a_-\lambda_-\tilde{\lambda}_-^{a_--1}+\lambda_-^{a_-}\big)
		=D^{(1)}(\mathbb{P}||\mathbb{Q}),\\
		B_{\Phi}(\tilde{\boldsymbol{\xi}},\boldsymbol{\xi})
		&=TC_+\Gamma(-a_+)\big((a_+-1)\lambda_+^{a_+}-a_+\tilde{\lambda}_+\lambda_+^{a_+-1}+\tilde{\lambda}_+^{a_+}\big)\nonumber\\
		&\quad+TC_-\Gamma(-a_-)\big((a_--1)\lambda_-^{a_-}-a_-\tilde{\lambda}_-\lambda_-^{a_--1}+\tilde{\lambda}_-^{a_-}\big)
		=D^{(-1)}(\mathbb{P}||\mathbb{Q}).
	\end{align}
	With the convention that $B_{\Phi}(\boldsymbol{\xi},\tilde{\boldsymbol{\xi}})$ is the canonical divergence in the $1$-affine coordinates, the two Bregman divergences above are the canonical and dual canonical divergences, respectively.

	The curvature tensor of the $\alpha$-connection is defined as
	\begin{align}
	\label{curvature_alpha}
		R^{(\alpha)l}_{\phantom{(\alpha)l}kij}
		&=\partial_i\Gamma^{(\alpha)l}_{jk}-\partial_j\Gamma^{(\alpha)l}_{ik}
		+\Gamma^{(\alpha)m}_{jk}\Gamma^{(\alpha)l}_{im}
		-\Gamma^{(\alpha)m}_{ik}\Gamma^{(\alpha)l}_{jm}.
	\end{align}
	The only nonzero connection components of the GTS geometry are $\Gamma^{(\alpha)\,\lambda_\pm}_{\lambda_\pm\lambda_\pm}$, and each depends only on the corresponding coordinate $\lambda_\pm$. The mixed derivatives and mixed connection components in Eq.~(\ref{curvature_alpha}) vanish. We obtain
	\begin{align}
		R^{(\alpha)}=0, \qquad \alpha\in\mathbb{R}.
	\end{align}
	For $\alpha=0$, the $0$-connection is the Levi-Civita connection. The Riemann curvature tensor of the Fisher metric vanishes, and its Ricci and scalar curvatures also vanish.

\subsection{CTS processes}
	It is straightforward to obtain CTS geometry from all the results we obtained above for GTS processes by applying the CTS condition of Eq.~(\ref{cts_condition}). 
	
	The $\alpha$-divergence between CTS processes is obtained as
	\begin{align}
	\label{a_div_cts}
		D^{(\alpha )}(\mathbb{P}||\mathbb{Q})=\left\{ 
		\begin{array}{ll}
		\frac{4}{1-\alpha^2}\Big\{1-\exp\big[-TC\Gamma(-a)\big(\big( \frac{1-\alpha}{2} \lambda_+^{a}+\frac{1+\alpha}{2} \tilde{\lambda}_+^{a}-( \frac{1-\alpha}{2} \lambda_++\frac{1+\alpha}{2} \tilde{\lambda}_+)^{a}\big)\\
		+\big( \frac{1-\alpha}{2} \lambda_-^{a}+\frac{1+\alpha}{2} \tilde{\lambda}_-^{a}-( \frac{1-\alpha}{2} \lambda_-+\frac{1+\alpha}{2} \tilde{\lambda}_-)^{a}\big)\big)\big]\Big\} & (\alpha \neq \pm 1)\\ 
		TC\Gamma(-a)\Big(\big((a-1)\lambda_+^{a}-a\tilde{\lambda}_+\lambda_+^{a-1}+\tilde{\lambda}_+^{a}\big)
		+\big((a-1)\lambda_-^{a}-a\tilde{\lambda}_-\lambda_-^{a-1}+\tilde{\lambda}_-^{a}\big)\Big) & (\alpha =-1)\\
		TC\Gamma(-a)\Big(\big((a-1)\tilde{\lambda}_+^{a}-a\lambda_+\tilde{\lambda}_+^{a-1}+\lambda_+^{a}\big)
		+\big((a-1)\tilde{\lambda}_-^{a}-a\lambda_-\tilde{\lambda}_-^{a-1}+\lambda_-^{a}\big)\Big) & (\alpha =1)
		\end{array}
		\right..
	\end{align}
	It is easily verified that the $\alpha=-1$ case in Eq.~(\ref{a_div_cts}) matches Eq.~(\ref{kld_cts}), which was originally obtained by Kim and Lee \cite{kim2007relative}.
	
	The Fisher information matrix of CTS processes can be obtained in three ways: plugging the CTS condition of Eq.~(\ref{cts_condition}) into the Fisher information matrix of GTS processes, Eq.~(\ref{metric_gts}); using Eq.~(\ref{metric_ts_lambda}) or Eq.~(\ref{metric_ts_t}); directly deriving it from the $\alpha$-divergence of Eq.~(\ref{a_div_cts}). These ways produce the same metric tensor of the CTS geometry as
	\begin{align}
	\label{metric_cts}
		g_{ij}=\Bigg(
		\begin{array}{cc}
			\frac{TC\Gamma(2-a)}{\lambda_+^{2-a}} &0 \\ 
			0 & \frac{TC\Gamma(2-a)}{\lambda_-^{2-a}}
		\end{array}\Bigg),
	\end{align}
	where $i,j$ run for $\lambda_+$ and $\lambda_-$. 
	
	Similar to the Fisher information matrix, the Levi-Civita connection of the CTS geometry is calculated from Eq.~(\ref{lc_conn}) or by imposing the CTS condition of Eq.~(\ref{cts_condition}) to Eq.~(\ref{lc_conn_gts_p}) and Eq.~(\ref{lc_conn_gts_m}). The non-trivial components of the Levi-Civita connection for the CTS geometry are found as
	\begin{align}
	\label{lc_conn_cts_p}
		\Gamma^{LC}_{\lambda_+\lambda_+,\lambda_+}&=-\frac{1}{2}\frac{TC\Gamma(3-a)}{\lambda_+^{3-a}},\\
	\label{lc_conn_cts_m}
		\Gamma^{LC}_{\lambda_-\lambda_-,\lambda_-}&=-\frac{1}{2}\frac{TC\Gamma(3-a)}{\lambda_-^{3-a}}.
	\end{align}
		
	By plugging the CTS condition of Eq.~(\ref{cts_condition}) to Eq.~(\ref{a_conn_gts_p}) and Eq.~(\ref{a_conn_gts_m}) or using Eq.~(\ref{levy_measure_cts}) and Eq.~(\ref{conn_ts}) in Theorem~\ref{thm_geo_ts}, the non-trivial components of the $\alpha$-connection for CTS processes are given by
	\begin{align}
	\label{a_conn_cts_p}
		\Gamma^{(\alpha)}_{\lambda_+\lambda_+,\lambda_+}&=-\frac{1-\alpha}{2}\frac{TC\Gamma(3-a)}{ \lambda_+^{3-a}},\\
	\label{a_conn_cts_m}
		\Gamma^{(\alpha)}_{\lambda_-\lambda_-,\lambda_-}&=-\frac{1-\alpha}{2}\frac{TC\Gamma(3-a)}{ \lambda_-^{3-a}}.
	\end{align}
	All other components vanish. Similar to the $\alpha$-connection of the GTS geometry, by plugging $\alpha=0$ to Eq.~(\ref{a_conn_cts_p}) and Eq.~(\ref{a_conn_cts_m}), we obtain the same Levi-Civita connection, Eq.~(\ref{lc_conn_cts_p}) and Eq.~(\ref{lc_conn_cts_m}), respectively. It is also possible to obtain the same result from Eq.~(\ref{conn_ts_lambda}) or Eq.~(\ref{conn_ts_t}).

	Eq.~(\ref{a_conn_cts_p}) and Eq.~(\ref{a_conn_cts_m}) vanish for $\alpha=1$. The coordinates $(\lambda_+,\lambda_-)$ are affine coordinates for the $1$-connection, and the CTS geometry is e-flat.

	A potential function for the Fisher metric is obtained as
	\begin{align}
	\label{potential_cts}
		\Phi_{\textrm{CTS}}(\boldsymbol{\xi})=TC\Gamma(-a)(\lambda_+^a+\lambda_-^a).
	\end{align}
	Using $\Gamma(2-a)=a(a-1)\Gamma(-a)$, the Hessian of Eq.~(\ref{potential_cts}) satisfies
	\begin{align}
		\partial_{\lambda_\pm}^2\Phi_{\textrm{CTS}}
		=TC\Gamma(2-a)\lambda_\pm^{a-2}
		=g_{\lambda_\pm\lambda_\pm},
	\end{align}
	and the mixed derivatives vanish. These are the components of the Fisher information matrix in Eq.~(\ref{metric_cts}).

	From the potential function, the dual coordinates are obtained as
	\begin{align}
	\label{dual_coord_cts}
		\eta_\pm=\partial_{\lambda_\pm}\Phi_{\textrm{CTS}}
		=TCa\Gamma(-a)\lambda_\pm^{a-1}.
	\end{align}
	The same calculation as for the GTS geometry gives $\nabla^{(-1)}d\eta_\pm=0$. The coordinates $(\eta_+,\eta_-)$ are affine coordinates for the $(-1)$-connection. The CTS statistical manifold is dually flat.

	Using Eq.~(\ref{potential_cts}), the Bregman divergences for CTS processes are given by
	\begin{align}
		B_{\Phi}(\boldsymbol{\xi},\tilde{\boldsymbol{\xi}})
		&=TC\Gamma(-a)\Big((a-1)\tilde{\lambda}_+^a-a\lambda_+\tilde{\lambda}_+^{a-1}+\lambda_+^a\nonumber\\
		&\qquad\qquad\quad +(a-1)\tilde{\lambda}_-^a-a\lambda_-\tilde{\lambda}_-^{a-1}+\lambda_-^a\Big)
		=D^{(1)}(\mathbb{P}||\mathbb{Q}),\\
		B_{\Phi}(\tilde{\boldsymbol{\xi}},\boldsymbol{\xi})
		&=TC\Gamma(-a)\Big((a-1)\lambda_+^a-a\tilde{\lambda}_+\lambda_+^{a-1}+\tilde{\lambda}_+^a\nonumber\\
		&\qquad\qquad\quad +(a-1)\lambda_-^a-a\tilde{\lambda}_-\lambda_-^{a-1}+\tilde{\lambda}_-^a\Big)
		=D^{(-1)}(\mathbb{P}||\mathbb{Q}).
	\end{align}
	These are the canonical and dual canonical divergences of the CTS geometry.

	Similar to the GTS geometry, $R^{(\alpha)}=0$ for every $\alpha\in\mathbb{R}$. The Riemann, Ricci, and scalar curvatures of the Fisher metric also vanish.

\subsection{RDTS processes}
	In a similar way, we can derive the $\alpha$-divergence between RDTS processes from Theorem~\ref{thm_div_ts}.
	
	\begin{crl}
		Let $(X_t, \mathbb{P})_{t\in[0,T]}$ and $(X_t, \mathbb{Q})_{t\in[0,T]}$ be RDTS processes with parameters $\boldsymbol{\xi}=(a_+, a_-, C_+, C_-,\lambda_+,\lambda_-,m)$ and $\tilde{\boldsymbol{\xi}}=(\tilde{a}_+, \tilde{a}_-, \tilde{C}_+, \tilde{C}_-,\tilde{\lambda}_+,\tilde{\lambda}_-,\tilde{m})$, respectively. Suppose $\mathbb{P}$ and $\mathbb{Q}$ are equivalent measures satisfying Eq.~(\ref{emm_rdts}). Then we have the $\alpha$-divergence between two RDTS processes as
	\begin{align}
	\label{a_div_rdts}
	D^{(\alpha )}(\mathbb{P}||\mathbb{Q})=\left\{
	\begin{array}{ll}
		\begin{aligned}
		\frac{4}{1-\alpha^2}\Bigg\{1-\exp\Bigg[&-2^{-1-\frac{a_+}{2}}TC_+\Gamma(-\frac{a_+}{2})\Big(\frac{1-\alpha}{2}\lambda_+^{\frac{a_+}{2}}+\frac{1+\alpha}{2}\tilde{\lambda}_+^{\frac{a_+}{2}}\\
		&\qquad-\big(\frac{1-\alpha}{2}\lambda_++\frac{1+\alpha}{2}\tilde{\lambda}_+\big)^{\frac{a_+}{2}}\Big)\\
		&-2^{-1-\frac{a_-}{2}}TC_-\Gamma(-\frac{a_-}{2})\Big(\frac{1-\alpha}{2}\lambda_-^{\frac{a_-}{2}}+\frac{1+\alpha}{2}\tilde{\lambda}_-^{\frac{a_-}{2}}\\
		&\qquad-\big(\frac{1-\alpha}{2}\lambda_-+\frac{1+\alpha}{2}\tilde{\lambda}_-\big)^{\frac{a_-}{2}}\Big)\Bigg]\Bigg\}
		\end{aligned} & (\alpha\neq\pm1)\\
			2^{-1-\frac{a_+}{2}}TC_+\Gamma(-\frac{a_+}{2})\big((\frac{a_+}{2}-1)\lambda_+^{\frac{a_+}{2}}-\frac{a_+}{2}\tilde{\lambda}_+\lambda_+^{\frac{a_+}{2}-1}+\tilde{\lambda}_+^{\frac{a_+}{2}}\big)\\
			+2^{-1-\frac{a_-}{2}}TC_-\Gamma(-\frac{a_-}{2})\big((\frac{a_-}{2}-1)\lambda_-^{\frac{a_-}{2}}-\frac{a_-}{2}\tilde{\lambda}_-\lambda_-^{\frac{a_-}{2}-1}+\tilde{\lambda}_-^{\frac{a_-}{2}}\big) & (\alpha=-1)\\
			2^{-1-\frac{a_+}{2}}TC_+\Gamma(-\frac{a_+}{2})\big((\frac{a_+}{2}-1)\tilde{\lambda}_+^{\frac{a_+}{2}}-\frac{a_+}{2}\lambda_+\tilde{\lambda}_+^{\frac{a_+}{2}-1}+\lambda_+^{\frac{a_+}{2}}\big)\\
			+2^{-1-\frac{a_-}{2}}TC_-\Gamma(-\frac{a_-}{2})\big((\frac{a_-}{2}-1)\tilde{\lambda}_-^{\frac{a_-}{2}}-\frac{a_-}{2}\lambda_-\tilde{\lambda}_-^{\frac{a_-}{2}-1}+\lambda_-^{\frac{a_-}{2}}\big) & (\alpha=1)
	\end{array}\right..
	\end{align}
	\end{crl}
	\begin{proof}
	Let us start with $\alpha\neq\pm1$. From Eq.~(\ref{hellinger_levy}) and Eq.~(\ref{a_div_measure_ts}) in Theorem~\ref{thm_div_ts}, the positive $x$ term of $-\log H_{\frac{1-\alpha}{2}}(\mathbb{P},\mathbb{Q})$ is given by
	\begin{align}
		&-\log H_{\frac{1-\alpha}{2}}(\mathbb{P},\mathbb{Q})|_{+}=T\int_0^\infty \Big(\frac{1-\alpha}{2}{\rm e}^{-\frac{(\lambda_+-\tilde{\lambda}_+)}{2}x^2} +\frac{1+\alpha}{2}- \Big({\rm e}^{-\frac{(\lambda_+-\tilde{\lambda}_+)}{2}x^2}\Big)^{\frac{1-\alpha}{2}}\Big)\frac{C_+{\rm e}^{-\frac{\tilde{\lambda}_+}{2} x^2}}{x^{a_++1}} dx\nonumber\\
	&=TC_+\int_0^\infty \Big(\frac{1-\alpha}{2}\frac{{\rm e}^{-\frac{\lambda_+}{2} x^2}}{x^{a_++1}} +\frac{1+\alpha}{2}\frac{{\rm e}^{-\frac{\tilde{\lambda}_+}{2} x^2}}{x^{a_++1}} -\frac{{\rm e}^{-\frac{1}{2}(\frac{1-\alpha}{2}\lambda_++\frac{1+\alpha}{2}\tilde{\lambda}_+) x^2}}{x^{a_++1}} \Big)dx\nonumber\\
	&=\frac{TC_+}{2}\int_0^\infty \Big(\frac{1-\alpha}{2}\frac{(\frac{\lambda_+}{2})^{a_+/2}{\rm e}^{-t}}{t^{a_+/2+1}} +\frac{1+\alpha}{2}\frac{(\frac{\tilde{\lambda}_+}{2})^{a_+/2}{\rm e}^{-t}}{t^{a_+/2+1}} -\frac{(\frac{1-\alpha}{2}\lambda_++\frac{1+\alpha}{2}\tilde{\lambda}_+)^{a_+/2}{\rm e}^{-t}}{t^{a_+/2+1}} \Big)dt\nonumber\\
	&=2^{-1-\frac{a_+}{2}}TC_+\Gamma(-\frac{a_+}{2})\Big(\frac{1-\alpha}{2}\lambda_+^{\frac{a_+}{2}}+\frac{1+\alpha}{2}\tilde{\lambda}_+^{\frac{a_+}{2}}-\big(\frac{1-\alpha}{2}\lambda_++\frac{1+\alpha}{2}\tilde{\lambda}_+\big)^{\frac{a_+}{2}}\Big).\nonumber
	\end{align}
	After a similar calculation is done for the negative $x$ part, we obtain
	\begin{align}
		-\log H_{\frac{1-\alpha}{2}}(\mathbb{P},\mathbb{Q})=&2^{-1-\frac{a_+}{2}}TC_+\Gamma(-\frac{a_+}{2})\Big(\frac{1-\alpha}{2}\lambda_+^{\frac{a_+}{2}}+\frac{1+\alpha}{2}\tilde{\lambda}_+^{\frac{a_+}{2}}-\big(\frac{1-\alpha}{2}\lambda_++\frac{1+\alpha}{2}\tilde{\lambda}_+\big)^{\frac{a_+}{2}}\Big)\nonumber\\
		&+2^{-1-\frac{a_-}{2}}TC_-\Gamma(-\frac{a_-}{2})\Big(\frac{1-\alpha}{2}\lambda_-^{\frac{a_-}{2}}+\frac{1+\alpha}{2}\tilde{\lambda}_-^{\frac{a_-}{2}}-\big(\frac{1-\alpha}{2}\lambda_-+\frac{1+\alpha}{2}\tilde{\lambda}_-\big)^{\frac{a_-}{2}}\Big).\nonumber
	\end{align}
	Substituting this expression into Eq.~(\ref{a_div_measure_ts}) gives the $\alpha$-divergence for $\alpha\neq\pm1$ in Eq.~(\ref{a_div_rdts}).

	For $\alpha=-1$, there are two ways to prove the result. The first way is to plug Eq.~(\ref{radon_nikodym_rdts}) to Eq.~(\ref{a_div_measure_ts}) (or Eq.~(\ref{log_radon_nikodym_rdts}) to Eq.~(\ref{a_div_measure_ts_psi})). The second way is to apply L'Hopital's rule to $\alpha \neq \pm1$ in the limit of $\alpha\to-1$. Both methods give the same result.
	
	We show the first method. Similar to $\alpha \neq \pm1$, we can start with the positive $x$ term of the Kullback--Leibler divergence:
	\begin{align}
		D^{(-1)}(\mathbb{P}||\mathbb{Q})|_{+}&=TC_+\int_0^\infty \Big(\frac{{\rm e}^{-\frac{\lambda_+}{2}x^2}}{x^{a_++1}}\big(-\frac{(\lambda_+-\tilde{\lambda}_+)}{2}x^2-1\big) + \frac{{\rm e}^{-\frac{\tilde{\lambda}_+}{2}x^2}}{x^{a_++1}}\Big)dx\nonumber\\
		&=2^{-1-\frac{a_+}{2}}TC_+\Gamma(-\frac{a_+}{2})\big(\frac{a_+}{2}\lambda_+^{\frac{a_+}{2}-1}(\lambda_+-\tilde{\lambda}_+)-\lambda_+^{\frac{a_+}{2}}+\tilde{\lambda}_+^{\frac{a_+}{2}}\big)\nonumber\\
		&=2^{-1-\frac{a_+}{2}}TC_+\Gamma(-\frac{a_+}{2})\big((\frac{a_+}{2}-1)\lambda_+^{\frac{a_+}{2}}-\frac{a_+}{2}\tilde{\lambda}_+\lambda_+^{\frac{a_+}{2}-1}+\tilde{\lambda}_+^{\frac{a_+}{2}}\big).\nonumber
	\end{align}
	
	After repeating the same step for negative $x$, the Kullback--Leibler divergence between two RDTS processes is given by
	\begin{align}
		D^{(-1)}(\mathbb{P}||\mathbb{Q})=&2^{-1-\frac{a_+}{2}}TC_+\Gamma(-\frac{a_+}{2})\big((\frac{a_+}{2}-1)\lambda_+^{\frac{a_+}{2}}-\frac{a_+}{2}\tilde{\lambda}_+\lambda_+^{\frac{a_+}{2}-1}+\tilde{\lambda}_+^{\frac{a_+}{2}}\big)\nonumber\\
		&+2^{-1-\frac{a_-}{2}}TC_-\Gamma(-\frac{a_-}{2})\big((\frac{a_-}{2}-1)\lambda_-^{\frac{a_-}{2}}-\frac{a_-}{2}\tilde{\lambda}_-\lambda_-^{\frac{a_-}{2}-1}+\tilde{\lambda}_-^{\frac{a_-}{2}}\big).\nonumber
	\end{align}
		
	For $\alpha=1$, we can use the same methods mentioned above. Additionally, we can use $\alpha$-duality with $\alpha=-1$. By following a similar calculation, the $\alpha$-divergence is obtained as
	\begin{align}
		D^{(1)}(\mathbb{P}||\mathbb{Q})=&D^{(-1)}(\mathbb{Q}||\mathbb{P})\nonumber\\
		=&2^{-1-\frac{a_+}{2}}TC_+\Gamma(-\frac{a_+}{2})\big((\frac{a_+}{2}-1)\tilde{\lambda}_+^{\frac{a_+}{2}}-\frac{a_+}{2}\lambda_+\tilde{\lambda}_+^{\frac{a_+}{2}-1}+\lambda_+^{\frac{a_+}{2}}\big)\nonumber\\
		&+2^{-1-\frac{a_-}{2}}TC_-\Gamma(-\frac{a_-}{2})\big((\frac{a_-}{2}-1)\tilde{\lambda}_-^{\frac{a_-}{2}}-\frac{a_-}{2}\lambda_-\tilde{\lambda}_-^{\frac{a_-}{2}-1}+\lambda_-^{\frac{a_-}{2}}\big)\nonumber.
	\end{align}
	\end{proof}
	
	It is noteworthy that the $\alpha$-divergence of RDTS processes is similar to those of CTS processes and GTS processes. It has additional scale factors, $2^{-1-\frac{a_\pm}{2}}$, and $a_{\pm}$ are replaced with $a_{\pm}/2$, respectively.

	We next consider the information geometry of RDTS processes.

	From Eq.~(\ref{levy_measure_rdts}) and Eq.~(\ref{metric_ts}) in Theorem~\ref{thm_geo_ts}, the Fisher information matrix of RDTS processes is found in the coordinate system of $\boldsymbol{\xi}=(\lambda_+,\lambda_-)$ as
	\begin{align}
	\label{metric_rdts}
		g_{ij}=\Bigg(
		\begin{array}{cc}
			\frac{2^{-1-\frac{a_+}{2}}TC_+\Gamma(2-\frac{a_+}{2})}{\lambda_+^{2-\frac{a_+}{2}}} &0 \\ 
			0 & \frac{2^{-1-\frac{a_-}{2}}TC_-\Gamma(2-\frac{a_-}{2})}{\lambda_-^{2-\frac{a_-}{2}}}
		\end{array}\Bigg),
	\end{align}
	where $i, j$ run for $\lambda_+$ and $\lambda_-$. We can obtain the same Fisher information matrix from Eq.~(\ref{metric_ts_lambda}) or Eq.~(\ref{metric_ts_t}).
	
	The diagonal metric tensor is explainable due to the same rationale regarding the indicator function with CTS processes and GTS processes. Additionally, based on the parametrization of the RDTS distribution, $2-\frac{a_{\pm}}{2}$ are positive.
	
	By Eq.~(\ref{lc_conn}), the non-trivial Levi-Civita connection components of the RDTS geometry are obtained as
	\begin{align}
	\label{lc_conn_rdts_p}
		\Gamma^{LC}_{\lambda_+\lambda_+,\lambda_+}&=-\frac{1}{2}\frac{2^{-1-\frac{a_+}{2}}TC_+\Gamma(3-\frac{a_+}{2})}{\lambda_+^{3-\frac{a_+}{2}}},\\
	\label{lc_conn_rdts_m}
		\Gamma^{LC}_{\lambda_-\lambda_-,\lambda_-}&=-\frac{1}{2}\frac{2^{-1-\frac{a_-}{2}}TC_-\Gamma(3-\frac{a_-}{2})}{\lambda_-^{3-\frac{a_-}{2}}}.
	\end{align}

	From Eq.~(\ref{levy_measure_rdts}) and Eq.~(\ref{conn_ts}) in Theorem~\ref{thm_geo_ts}, the non-trivial $\alpha$-connection components of the RDTS manifolds are given by
	\begin{align}
	\label{a_conn_rdts_p}
		\Gamma^{(\alpha)}_{\lambda_+\lambda_+,\lambda_+}&=-\frac{1-\alpha}{2}\frac{2^{-1-\frac{a_+}{2}}TC_+\Gamma(3-\frac{a_+}{2})}{ \lambda_+^{3-\frac{a_+}{2}}},\\
	\label{a_conn_rdts_m}
		\Gamma^{(\alpha)}_{\lambda_-\lambda_-,\lambda_-}&=-\frac{1-\alpha}{2}\frac{2^{-1-\frac{a_-}{2}}TC_-\Gamma(3-\frac{a_-}{2})}{ \lambda_-^{3-\frac{a_-}{2}}}.
	\end{align}
	All other components vanish. By plugging $\alpha=0$ to Eq.~(\ref{a_conn_rdts_p}) and Eq.~(\ref{a_conn_rdts_m}), we can obtain the same Levi-Civita connection, Eq.~(\ref{lc_conn_rdts_p}) and Eq.~(\ref{lc_conn_rdts_m}). The same $\alpha$-connection can be calculated from Eq.~(\ref{conn_ts_lambda}) or Eq.~(\ref{conn_ts_t}).

	For $\alpha=1$, Eq.~(\ref{a_conn_rdts_p}) and Eq.~(\ref{a_conn_rdts_m}) vanish. The coordinates $(\lambda_+,\lambda_-)$ are affine coordinates for the $1$-connection, and the RDTS geometry is e-flat.

	A potential function for the Fisher metric is obtained as
	\begin{align}
	\label{potential_rdts}
		\Phi_{\textrm{RDTS}}(\boldsymbol{\xi})
		&=2^{-1-\frac{a_+}{2}}TC_+\Gamma(-\frac{a_+}{2})\lambda_+^{\frac{a_+}{2}}
		+2^{-1-\frac{a_-}{2}}TC_-\Gamma(-\frac{a_-}{2})\lambda_-^{\frac{a_-}{2}}.
	\end{align}
	Using $\Gamma(2-a/2)=(a/2)(a/2-1)\Gamma(-a/2)$, the Hessian of Eq.~(\ref{potential_rdts}) satisfies
	\begin{align}
		\partial_{\lambda_\pm}^2\Phi_{\textrm{RDTS}}
		=2^{-1-\frac{a_\pm}{2}}TC_\pm\Gamma(2-\frac{a_\pm}{2})\lambda_\pm^{\frac{a_\pm}{2}-2}
		=g_{\lambda_\pm\lambda_\pm},
	\end{align}
	and the mixed derivatives vanish. These are the components of the Fisher information matrix in Eq.~(\ref{metric_rdts}).

	From the potential function, the dual coordinates are obtained as
	\begin{align}
	\label{dual_coord_rdts}
		\eta_\pm=\partial_{\lambda_\pm}\Phi_{\textrm{RDTS}}
		=2^{-1-\frac{a_\pm}{2}}TC_\pm\frac{a_\pm}{2}\Gamma(-\frac{a_\pm}{2})\lambda_\pm^{\frac{a_\pm}{2}-1}.
	\end{align}
	The same calculation as for the GTS geometry gives $\nabla^{(-1)}d\eta_\pm=0$. The coordinates $(\eta_+,\eta_-)$ are affine coordinates for the $(-1)$-connection. The RDTS statistical manifold is dually flat.

	Using Eq.~(\ref{potential_rdts}), the Bregman divergences for RDTS processes are given by
	\begin{align}
		B_{\Phi}(\boldsymbol{\xi},\tilde{\boldsymbol{\xi}})
		&=2^{-1-\frac{a_+}{2}}TC_+\Gamma(-\frac{a_+}{2})\Big((\frac{a_+}{2}-1)\tilde{\lambda}_+^{\frac{a_+}{2}}-\frac{a_+}{2}\lambda_+\tilde{\lambda}_+^{\frac{a_+}{2}-1}+\lambda_+^{\frac{a_+}{2}}\Big)\nonumber\\
		&\quad+2^{-1-\frac{a_-}{2}}TC_-\Gamma(-\frac{a_-}{2})\Big((\frac{a_-}{2}-1)\tilde{\lambda}_-^{\frac{a_-}{2}}-\frac{a_-}{2}\lambda_-\tilde{\lambda}_-^{\frac{a_-}{2}-1}+\lambda_-^{\frac{a_-}{2}}\Big)
		=D^{(1)}(\mathbb{P}||\mathbb{Q}),\\
		B_{\Phi}(\tilde{\boldsymbol{\xi}},\boldsymbol{\xi})
		&=2^{-1-\frac{a_+}{2}}TC_+\Gamma(-\frac{a_+}{2})\Big((\frac{a_+}{2}-1)\lambda_+^{\frac{a_+}{2}}-\frac{a_+}{2}\tilde{\lambda}_+\lambda_+^{\frac{a_+}{2}-1}+\tilde{\lambda}_+^{\frac{a_+}{2}}\Big)\nonumber\\
		&\quad+2^{-1-\frac{a_-}{2}}TC_-\Gamma(-\frac{a_-}{2})\Big((\frac{a_-}{2}-1)\lambda_-^{\frac{a_-}{2}}-\frac{a_-}{2}\tilde{\lambda}_-\lambda_-^{\frac{a_-}{2}-1}+\tilde{\lambda}_-^{\frac{a_-}{2}}\Big)
		=D^{(-1)}(\mathbb{P}||\mathbb{Q}).
	\end{align}
	The two Bregman divergences give the canonical and dual canonical divergences of the RDTS geometry.

	The same diagonal structure as for the GTS geometry gives $R^{(\alpha)}=0$ for every $\alpha\in\mathbb{R}$. The Riemann, Ricci, and scalar curvatures of the Fisher metric also vanish.

\section{Statistical applications}
\label{sec_apps}
	In this section, we introduce statistical applications of information geometry for GTS processes. Since the geometry of CTS processes and RDTS processes is similar to the GTS geometry, we can apply the findings regarding the GTS processes in this section to those two processes.

\subsection{Bias reduction}
	First of all, we apply our findings to bias reduction in maximum likelihood estimates. Firth \cite{firth1993bias} suggested a more systematic bias reduction in maximum likelihood estimation. For canonical exponential families, the penalized log-likelihood given in Firth \cite{firth1993bias} can be written as
	\begin{align}
	\label{penalized_ll}
		l^{*}(\boldsymbol{\xi})=l(\boldsymbol{\xi})+\log \mathcal{J},
	\end{align}
	where $l(\boldsymbol{\xi})$ is the log-likelihood and $\mathcal{J}$ is the Jeffreys prior. From this penalized log-likelihood, we can estimate bias-reduced parameters.
	
	The Jeffreys prior, a non-informative prior distribution, is defined as follows \cite{jeffreys1946invariant}:
	\begin{align}
		\mathcal{J}(\boldsymbol{\xi})\propto g^{\frac{1}{2}},
	\end{align}
	where $g$ is the determinant of the Fisher information matrix which is the metric tensor of information geometry.
	
	From the Fisher information matrix of GTS processes, Eq.~(\ref{metric_gts}), it is easy to obtain the Jeffreys prior of GTS processes as
	\begin{align}
		\label{jeffreys_gts}
		\mathcal{J}(\boldsymbol{\xi})\propto T\sqrt{\frac{C_+\Gamma(2-a_+)}{\lambda_+^{2-a_+}} \frac{C_-\Gamma(2-a_-)}{\lambda_-^{2-a_-}}}.
	\end{align}
	
	The Jeffreys prior in Eq.~(\ref{jeffreys_gts}) can be used in the penalized likelihood of Eq.~(\ref{penalized_ll}) for a given likelihood model.

	For a numerical example, we test the Firth estimation for CTS processes. Let us consider a CTS process with $a=0.7$, $C=0.1$, $T=1$, and the true parameters $(\lambda_+,\lambda_-)=(1,1.5)$. We generate $n$ independent unit-horizon full-path observations and estimate $\lambda_+$ and $\lambda_-$ while the other parameters are fixed at their true values. The positive and negative jump components are independent, so we calculate the two one-sided likelihoods separately.

	For either one-sided component, let $\tilde{\lambda}$ be a reference parameter and let $x>0$ denote the jump magnitude. From Eq.~(\ref{log_radon_nikodym_gts}), $\psi(x)=-(\lambda-\tilde{\lambda})x$. Since $a=0.7<1$, Eq.~(\ref{rn_path}) can be applied to each of the $n$ independent paths. Let $N_i$ be the jump measure of $X^{(i)}$ and $\nu_{\tilde{\lambda}}(dx)=C{\rm e}^{-\tilde{\lambda}x}x^{-a-1}dx$. Then
	\begin{align*}
		\log\frac{d\mathbb{P}_{\lambda}^{\otimes n}}{d\mathbb{P}_{\tilde{\lambda}}^{\otimes n}}
		&=\sum_{i=1}^n\log\frac{d\mathbb{P}_{\lambda}}{d\mathbb{P}_{\tilde{\lambda}}}(X^{(i)})\\
		&=\sum_{i=1}^n\left\{\int_0^T\int\psi(x)N_i(dt,dx)-T\int\left({\rm e}^{\psi(x)}-1\right)\nu_{\tilde{\lambda}}(dx)\right\}\\
		&=-(\lambda-\tilde{\lambda})\sum_{i=1}^n\sum_{t\leq T}|\Delta X_t^{(i)}|-nTC\Gamma(-a)\left(\lambda^a-\tilde{\lambda}^a\right)\\
		&=-(\lambda-\tilde{\lambda})S-nTC\Gamma(-a)\left(\lambda^a-\tilde{\lambda}^a\right),
	\end{align*}
	where $S=\sum_{i=1}^n\sum_{t\leq T}|\Delta X_t^{(i)}|$ is the total jump variation of the corresponding one-sided component over the $n$ paths. With $\kappa=-TC\Gamma(-a)$ and the terms independent of $\lambda$ omitted, the log-likelihood is given by
	\begin{align*}
		l(\lambda)=-\lambda S+n\kappa\lambda^a+\mathrm{const}.
	\end{align*}
	The ordinary maximum likelihood estimator solves $-S+n\kappa a\lambda^{a-1}=0$. From Eq.~(\ref{jeffreys_gts}), $\partial_\lambda\log\mathcal{J}=-(2-a)/(2\lambda)$, and the interior Firth adjusted-score estimator associated with Eq.~(\ref{penalized_ll}) solves
	\begin{align}
		-S+n\kappa a\lambda^{a-1}-\frac{2-a}{2\lambda}=0.
	\end{align}
	For $n=50,100,200,500$, we repeat the simulation $5,000$ times and calculate the relative bias $(\mathbb{E}[\hat{\lambda}]-\lambda)/\lambda$ and the root mean squared error (RMSE). 
	
\begin{figure}[t]
		\centering
		\subfigure[Relative bias]{\includegraphics[width=0.47\textwidth]{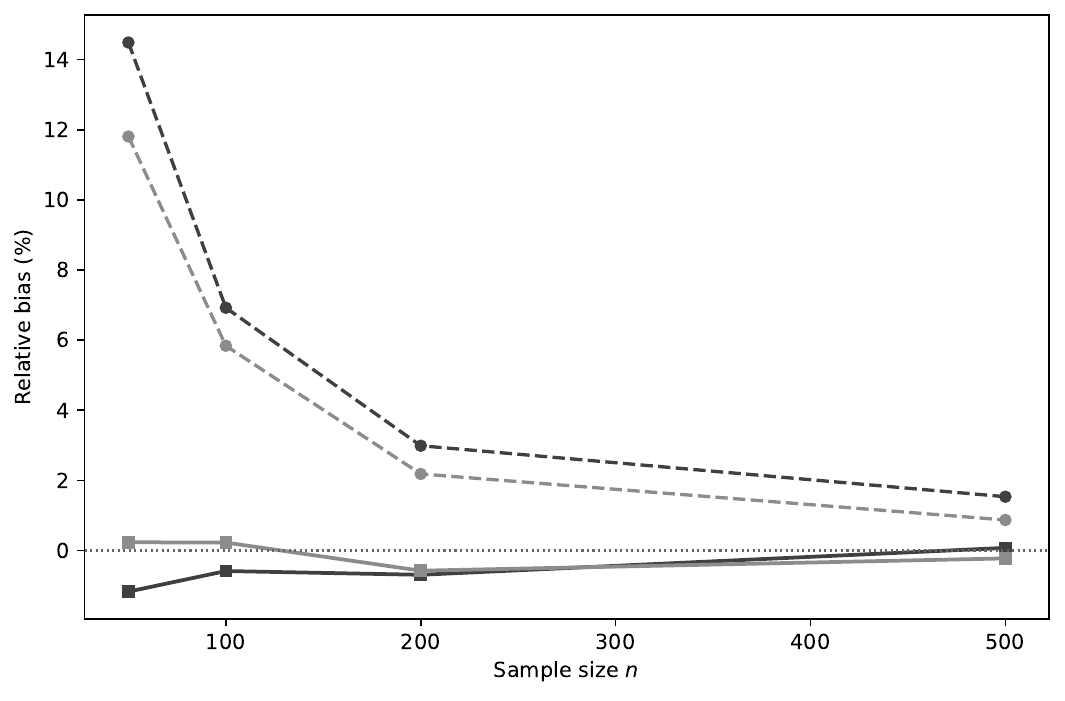}}\hfill
		\subfigure[RMSE]{\includegraphics[width=0.47\textwidth]{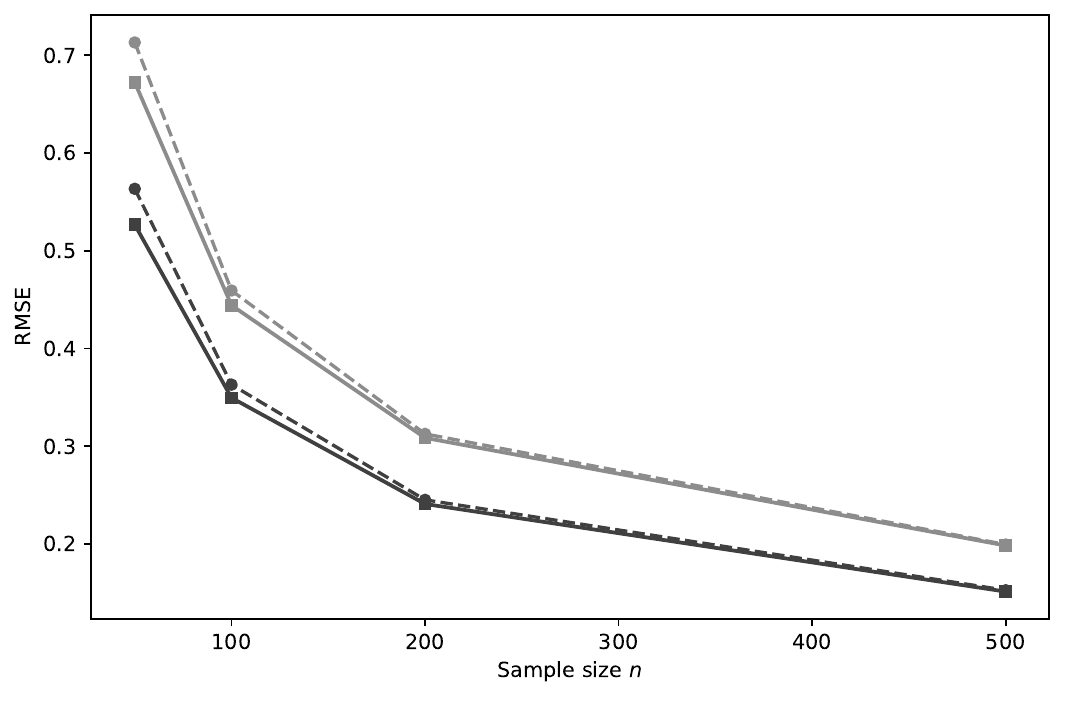}}
		\caption{Relative biases and RMSEs of the maximum likelihood and Firth estimators for $\lambda_+$ and $\lambda_-$ in the CTS model. Solid-square and dashed-circle lines denote the Firth and maximum likelihood estimators, respectively. Dark and light gray lines correspond to $\lambda_+$ and $\lambda_-$, respectively. The results are based on $5,000$ Monte Carlo replications with $a=0.7$, $C=0.1$, $T=1$, and $(\lambda_+,\lambda_-)=(1,1.5)$.}
		\label{fig_firth_cts}
	\end{figure}

	The relative biases and RMSEs are shown in Figure~\ref{fig_firth_cts}. The MLEs have positive relative biases for small $n$. The relative biases of the Firth estimators remain close to zero. Additionally, the Firth estimators have lower RMSEs than the MLEs for both $\lambda_+$ and $\lambda_-$ at all sample sizes. The differences in the RMSEs are larger for small $n$ and become smaller as $n$ increases.

\subsection{Bayesian predictive priors}
	Another application is finding Bayesian predictive priors for GTS processes. From the metric tensor of the GTS geometry, we calculate several ansatze satisfying the superharmonicity condition for obtaining shrinkage priors \cite{komaki2006shrinkage}.
		
	Komaki \cite{komaki2006shrinkage} derived Bayesian predictive priors outperforming the Jeffreys prior. The Bayesian predictive prior is represented with
	\begin{align}
		\tilde{\mathcal{J}}=\phi^2 \mathcal{J},
	\end{align}
	where $\phi$ is a positive superharmonic function and $\mathcal{J}$ is the Jeffreys prior.
	
	Komaki-style shrinkage priors for time series models and signal filters are already well-known. We have explicit forms of the predictive priors for AR models \cite{komaki2006shrinkage,tanaka2018superharmonic}, ARMA models \cite{choi2015kahlerian, choi2015geometric, oda2021shrinkage}, and signal processing filters \cite{choi2015geometric}.
	
	Similar to the previous work on finding ansatze for Bayesian predictive priors of signal processing filters \cite{choi2015geometric}, we can calculate ansatze for GTS processes. Let us try to find ansatze for superharmonic functions $\phi$ on the GTS statistical manifold. 
	
	The first ansatz is a function in $\lambda_+$:
	\begin{align}
		\phi_1=\lambda_+^k,
	\end{align}
	where $0<k<a_+/2$. It is easy to verify that it is a positive superharmonic function.

	Similarly, we can construct another ansatz in $\lambda_-$:	
	\begin{align}
		\phi_2=\lambda_-^l,
	\end{align}
	where $0<l<a_-/2$. It is also straightforward to check the superharmonicity of the ansatz.
	
	We can consider a positive linear combination and the product of $\phi_1$ and $\phi_2$:
	\begin{align}
		\phi_3&= c_1\phi_1+c_2\phi_2,\\
		\phi_4&= \phi_1\phi_2,
	\end{align}
	where $c_1>0$ and $c_2>0$. It is straightforward to check that these are also positive superharmonic functions.
		
\section{Conclusion}
	We derived the information geometry of tempered stable processes. It is noteworthy that finding the $\alpha$-divergence between two tempered stable processes is a two-fold generalization of the work by Kim and Lee \cite{kim2007relative}: from CTS processes to tempered stable processes and from the Kullback--Leibler divergence to the $\alpha$-divergence. 
	
	From the $\alpha$-divergence of tempered stable processes, we obtained the information geometry of tempered stable processes: Fisher information matrix and $\alpha$-connection. Since the GTS process parameters for the upper and lower tails are not entangled with each other, the metric tensor of the GTS process geometry is diagonal. Additionally, we further obtained dual flatness, potential functions, canonical divergences, and vanishing curvature for the GTS, CTS, and RDTS geometries. 

	With these geometric objects, we also provided several information-geometric applications for GTS processes. First of all, we considered bias reduction in maximum likelihood estimation for parameter estimation. The CTS example also shows the bias reduction numerically. Additionally, candidate shrinkage priors based on positive superharmonic functions were obtained. 

\section*{Acknowledgment}
	We thank Young Shin Kim for useful discussions on tempered stable processes and distributions and the relative entropy of CTS processes.
	
\bibliographystyle{plain}
\bibliography{ig_gts}

\end{document}